\documentclass[11pt]{amsart}
\usepackage{amsmath,amsfonts,amsthm}
\usepackage{amssymb,latexsym}
\usepackage{rotating}
\usepackage[colorlinks=true, linkcolor=blue, citecolor=red]{hyperref}
\usepackage{graphicx}
\usepackage[all]{xy}
\usepackage{tikz}
\usepackage{layout}
\theoremstyle{plain}
\newtheorem{theorem}{Theorem}[section]
\newtheorem{lemma}[theorem]{Lemma}
\newtheorem{proposition}[theorem]{Proposition}
\newtheorem{corollary}[theorem]{Corollary}

\theoremstyle{definition}

\newtheorem{definition}[theorem]{Definition}
\newtheorem{example}[theorem]{Example}

\newtheorem*{proofb}{Proof of Theorem~\ref{thm_main}}
\newtheorem*{proofc}{Proof of Theorem~\ref{thm_dim_k0_iso}}
\newtheorem*{proofd}{Proof of Theorem~\ref{thm_k0_dim_iso}}
\newtheorem*{proofe}{Proof of Corollary~\ref{cor_mu_system}}

\newtheorem*{notation}{Notation}
\newtheorem{remark}[theorem]{Remark}

\numberwithin{equation}{section}

\newtheorem{defx}{Definition}

\newtheorem{thmx}[defx]{Theorem}

\newtheorem{corx}[defx]{Corollary}

\newenvironment{proof_thmxx}{
\begin{proofc}}
{$\square$
\end{proofc}}

\newenvironment{proof_thmxxx}{
\begin{proofd}}
{$\square$
\end{proofd}}

\newenvironment{proof_corx}{
\begin{proofe}}
{$\square$
\end{proofe}}

\RequirePackage[normalem]{ulem} 
\RequirePackage{color}
\definecolor{RED}{rgb}{1,0,0}
\definecolor{BLUE}{rgb}{0,0,1}

\title{
Maximal UHF~subalgebras of certain C*-algebras   
}

\author{Nasser Golestani and Saeid Maleki Ouche}

\address[\textbf{Nasser Golestani}]{Department of Pure Mathematics, 
Faculty of Mathematical Sciences, Tarbiat Modares University, Tehran \\ 
Iran}
\email{n.golestani@modares.ac.ir}

\address[\textbf{Saeid Maleki Ouche}]{Department of Pure Mathematics, 
Faculty of Mathematical Sciences, Tarbiat Modares University, Tehran \\ 
Iran} 
\email{saeidmalekiouche@modares.ac.ir}

\begin{document}
\maketitle
\begin{abstract}
A well-known result in dynamical systems asserts that 
any Cantor minimal system \( (X, T) \) has 
a maximal rational equicontinuous factor \( (Y, S) \) 
which is in fact an odometer, and realizes the rational subgroup 
of the \( K_{0} \)-group of \( (X, T) \), that is,
\( \mathbb{Q}(K^{0}(X, T), 1) \cong K^{0}(Y, S) \). We introduce 
the notion of a maximal UHF~subalgebra and use it to obtain 
the C*-algebraic anolog of this result.
We say a UHF~subalgebra ${B}$ of a unital C*-algebra $A$
is a maximal UHF~subalgebra if it contains the unit of $A$ 
and any other such C*-subalgebra of $A$ embeds unitaly into 
$B$. We prove that if $K_{0}(A)$ is unperforated and has 
a certain $K_{0}$-lifting property, then ${B}$ exists and is 
unique up to isomorphism, in particular, all 
simple separable unital C*-algebras with   
tracial rank zero and all unital Kirchberg algebras whose 
$K_{0}$-groups are unperforated, have 
a maximal UHF~subalgebra. Not every unital C*-algebra has 
a maximal UHF~subalgebra, for instance, 
the unital universal free product 
$\mathrm{M}_{2} \ast_{r} \mathrm{M}_{3}$. 
As an application, we give a C*-algebraic realization of 
the rational subgroup $\mathbb{Q}(G, u)$ of 
any dimension group $G$ with order unit $u$, that is, there is
a simple unital AF~algebra (and a unital Kirchberg algebra) $A$ 
with a maximal UHF~sublgebra $B$ such that 
$(G, u) \cong \left(K_{0}(A), [1]_{0}\right)$ and 
$\mathbb{Q}(G, u) \cong K_{0}(B)$.
\end{abstract}

\tableofcontents
\section{Introduction}\label{sec_intro}

In operator algebras, certain subalgebras play an important role. 
For instance, Cartan subalgebras of von Neumann algebras and 
C*-algebras \cite{re08}, and large subalgebras of 
simple unital C*-algebras \cite{ph14}. In this paper, we consider 
maximal UHF~subalgebras of unital C*-algebras. Our first motivation 
is to give a C*-algebraic realization of the rational subgroup of 
a dimension group in such a way that it has a suitable relation to 
the dynamical realization (see Theorems~\ref{thm_dim_k0_iso} 
and \ref{cor_mu_system}). Moreover, as UHF~algebras are 
well understand in operators algebras a maximal UHF~subalgebra of 
a unital C*-algebra \( A \) (if exists), may be useful to understand 
some aspects of the structure of \( A \), in particular its K-theory 
(see Theorem~\ref{thm_k0_dim_iso}).

Dimension groups were introduced by G. A. Elliott for 
the classification of AF~algebras \cite{el76}. Since then they 
became a powerful tool to study the K-theory of both 
C*-algebras and Cantor minimal systems \cite{ef81, gps95}.  
The \emph{rational subgroup} of a dimension group $G$ with 
order unit $u$ \cite{or97, ghh17} is defined by
\[ \mathbb{Q}(G, u) = \left\{ g \in G : mg = qu \ 
\text{for some \( m \in \mathbb{N} \) and} \ q \in \mathbb{Z}
\ \right\}. \]

The dynamical realization of the rational subgroup was given 
using the maximal rational equicontinuous factor of 
a Cantor minimal system. More precisely, for every 
dimension group $G$ with order unit $u$ there is 
a Cantor minimal system $(X, T)$ such that 
$G \cong K^{0}(X, T)$ and 
$\mathbb{Q}(G, u) \cong K^{0}(Y, S)$ 
where $(Y, S)$ is an odometer and is
the maximal rational equicontinuous factor of $(X, T)$ 
\cite{gps95, ho16}.

\medskip

Our first aim is to find a suitable C*-algebraic realization of 
the rational subgroup of an ordered Abelian group. 
As odometers corresponds to UHF~algebras (since both have 
Bratteli diagrams with one vector at each level \cite{aeg21, br72, rll00}) 
and dynamical factors corresponds to C*-subalgebras,
we introduce the following notion.

\begin{defx}\label{def_mu}
A UHF subalgebra $B$ of a unital C*-algebra $A$ is
a \emph{maximal UHF~subalgebra} of $A$ if $1_{B} = 1_{A}$ 
and for any UHF C*-subalgebra $D$ of $A$ with 
$1_{D} = 1_{A}$, there exists a unital embedding from $D$ to 
$B.$ If such a $B$ exists (which is necessarily unique 
up to isomorphism), we denote it by $MU(A)$.
\end{defx}

\medskip
 
A maximal UHF~subalgebra of the following C*-algebras is 
isomorphic to $\mathbb{C}$: the Jiang-Su algebra 
$\mathcal{Z}$, the Toeplitz algebra $\mathcal{T}$, $C(X)$ for 
any compact Hausdorff space $X$, and 
the unitization $\tilde{A}$ of any nonunital C*-algebra $A$. 
On the other hand, a maximal UHF~subalgebra of 
the Cuntz algebra $\mathcal{O}_{2}, B(\mathcal{H})$, and 
the Calkin algebra $\mathcal{Q}(\mathcal{H})$ for any 
infinite dimensional Hilbert space $\mathcal{H}$ is 
the universal UHF~algebra $\mathcal{Q}$. See 
Proposition~\ref{prop_mu} for a list of examples.

\medskip

In the following theorem, we determine a class of C*-algebras 
having a maximal UHF~subalgebra.
We say that a unital C*-algebra $A$ has 
the \emph{$K_{0}$-lifting property for UHF~algebras} if 
the existence of an injective positive order unit preserving 
group homomorphism $K_{0}(D) \to K_{0}(A)$ where $D$ is 
a UHF~algebra, implies the existence of 
a unital $*$-homomorphism $D \to A$. 

\begin{thmx}\label{thm_main}
Every unital C*-algebra with $K_{0}$-lifting property for 
UHF~algebras whose $K_{0}$-group is unperforated, has  
a maximal UHF~subalgebra.  
\end{thmx}

For instance, all simple separable unital C*-algebras with 
tracial rank zero and all unital Kirchberg algebras whose 
$K_{0}$-groups are unperforated, have a maximal UHF~subalgebra. 
In Section~\ref{sec_mu_bra}, we give a combinatorial method 
based on Bratteli diagrams to construct a maximal UHF~subalgebra 
for any unital AF~algebra.

To prove Theorem~\ref{thm_main}, first we introduce 
in subsection~\ref{subsec_propd} the notion of Property~(D) 
for an ordered Abelian group with a distinguished order unit
\( (G, G^{+}, u) \) which says that if \( m \vert u \) and \( n | u \) 
for co-prime natural numbers \( m \) and \( n \), then \( mn | u \). 
Every weakly unperforated ordered Abelian group has this property. 
Next, we show that this property guarantees the existence of 
the largest supernatural number \( N = N(G, u) \) dividing \( u \) 
(Theorem~\ref{thm_max_propd}). Then we obtain an embedding 
\( Q(N) \to G \), when \( G \) is unperforated. Finally, if \( A \) has 
the \( K_{0} \)-lifting property for UHF~algebras and \( K_{0}(A) \) is 
unperforated, we take \( G = K_{0}(A) \) and we complete the proof 
in subsection~\ref{subsec_k0lift} 

If $A$ and $B$ are unital C*-algebras such that $A$ has 
a maximal UHF~subalgebra and there are 
unital $*$-homomorphisms $A \to B$ and $B \to A$, then 
$B$ has a maximal UHF~subalgebra which is isomorphic to that 
of $A$ (Proposition~\ref{prop_mu_homo}). In particular, this is 
the case if $A$ and $B$ are homotopy equivalent.

\begin{thmx}\label{thm_dim_k0_iso}
If $G$ is a dimension group with an order unit $u$, then there is
a unital AF~algebra (and a unital Kirchberg algebra) $A$ 
with a maximal UHF~sublgebra $B$ such that 
$(K_{0}(A), [1]_{0}) \cong (G, u)$ and 
$\mathbb{Q}(G, u) \cong K_{0}(B)$.
\end{thmx}

In fact, there is an uncountable family of 
pairwise nonisomorphic C*-algebras $A$ satisfying 
the preceding theorem. We can arrange this family to consist of 
simple unital AF~algebras or unital Kirchberg algebras 
(Theorems~\ref{thm_mu_af} and \ref{thm_mu_kirchberg}). 

\medskip

Our first application of these results is a C*-algebraic realization 
of the rational subgroup $\mathbb{Q}(G, u)$ of 
a dimension group $G$ with order unit $u$.

\begin{thmx}\label{thm_k0_dim_iso}
Let $A$ be a unital C*-algebra having 
a maximal UHF~subalgebra. If $(K_{0}(A), [1]_{0})$ is 
a dimension group then 
$K_{0}(MU(A)) \cong \mathbb{Q}(K_{0}(A), [1]_{0})$ 
as dimension groups with order unit. In particular, this is 
the case if $A$ is a unital AF~algebra.
\end{thmx}

The proof of these two results requires some ingredients:
the part of the Elliott classification program dealing with 
the range of the Elliott invariant, the isomorphism 
\( K_{0}(MU(A)) \cong Q(N(K_{0}(A), [1]_{0})) \)
already provided in the proof of Theorem~\ref{thm_main},
and a realization of the rational subgroup of a dimension group
\( (G, u) \) by \( Q(N(G, u), 1) \) given in 
Theorem~\ref{thm_grp_iso_2}.

As another application of these results, we are able to make
a connection between dynamical and C*-algebraic realizations 
of the rational subgroups of dimension groups as follows.

\begin{corx}\label{cor_mu_system}
Let $(X, T)$ be a Cantor minimal system with 
the maximal rational equicontinuous factor $(Y, S)$. Then 
$K^{0}(Y, S) \cong K_{0}(B)$ 
as dimension groups with order unit where $B$ is 
a maximal UHF~subalgebra of the C*-algebra crossed product 
$C(X) \rtimes_{T} \mathbb{Z}$.
\end{corx}

The structure of this paper is as follows. 
In Section~\ref{sec_mu_grp} we give some preliminaries 
on ordered Abelian groups, introduce Property~(D), and 
prove Theorem~\ref{thm_main}. Section~\ref{sec_exa} 
is devoted to the permanence properties and various examples 
of C*-algebras having maximal UHF~subalgebras. 
In Section~\ref{sec_rat_grp}, we prove 
Theorems~\ref{thm_dim_k0_iso}, \ref{thm_k0_dim_iso}, and 
\ref{cor_mu_system}. In the final section, we use Bratteli diagrams 
to give a constructive and combinatorial method to obtain 
a maximal UHF~subalgebra of a unital AF~algebra.

\pagestyle{plain} 

\section{Maximal UHF~subalgebras, Ordered Groups Approach}
\label{sec_mu_grp}

\begin{notation}
We use the following notation throughout this paper.
\hfill
\begin{enumerate}
\item 
$A^{+}$ denotes the unitization of a C*-algebra $A$ (adding 
a new identity even if $A$ is unital), while $A^{\sim} = A$ if 
$A$ is unital and $A^{\sim} = A^{+}$ if $A$ is nonunital.

\item
$\mathcal{K} = \mathcal{K}(\ell^{2})$ and 
$\mathrm{M}_{n} = M_{n}(\mathbb{C})$.

\item
We denote the universal UHF~algebra associated with 
the supernatural number $N = \left\{ \infty, \infty, \ldots \right\}$ 
by $\mathcal{Q}$.
	
\item
We write $A\sim_h B$ if $A$ and $B$ are homotopy equivalent 
C*-algebras.	

\item
For separable C*-algebras $A, B$, two $*$-homomorphisms
$\varphi ,\psi : A \to B$ are called approximately unitarily (a.u.) 
equivalent, denoted by $\varphi \approx_{a.u.} \psi$, if there is 
a sequence $(u_n)_{n = 1}^{\infty}$ of unitaries in $B^{\sim}$ 
such that 
$\lim_{n \to \infty} \Vert u^*_n \varphi(a) u_n - \psi (a) \Vert = 0$ 
for all $a \in A$.
\end{enumerate}
\end{notation}

\medskip

\subsection{Ordered Abelian Groups}\label{subsec_ord_grp}
In this subsection, we recall notions about 
ordered Abelian groups and UHF~algebras \cite{rll00}. A pair 
$(G, G^{+})$ is called an \emph{ordered Abelian group} if $G$
is an Abelian group, $G^{+}$ is a subset of $G$, and 
\[ G^{+} + G^{+} \subseteq G^{+}, \hspace{0.5 cm} 
G^{+} \cap (- G^{+}) = \{ 0 \}, \hspace{0.5 cm} 
G^{+} - G^{+} = G. \] 
The a relation $\leq$ on $G$ is defined by $x \leq y$ if 
$y - x \in G^{+}$. Note that some authors do not assume 
the third property above when defining 
an ordered Abelian group \cite[Page 82]{dp22}.

\medskip

An element $u$ in $G^{+}$ in an ordered Abelian group 
$(G, G^{+})$ is called an \emph{order unit} if for every $g$ 
in $G$ there is a positive integer $n$ with $-nu \leq g \leq nu$.
A triple $(G, G^{+}, u)$, where $(G, G^{+})$ is 
an ordered Abelian group and $u$ is an order unit, is called 
an \emph{ordered Abelian group with a distinguished order unit}.

Let $(G, G^{+})$ be an ordered Abelian group. If $x$ in $G$ 
for which $nx > 0$ for some $n \in \mathbb{N}$ satisfies 
$x > 0$, then $G$ is said to be \emph{weakly unperforated}.
Similarly, if $nx \geq 0$ implies $x \geq 0$, then $G$ is called 
\emph{unperforated}.
\medskip


Unless specified explicitly, we equip the ordered Abelian group 
$\mathbb{Z}^{d}$ with the natural cone
$(\mathbb{Z}^{+})^{d}$ where 
$\mathbb{Z}^{+} = \left\{ 0, 1, 2, \ldots \right\}$, and with order unit
$\left(1, 1, \ldots, 1\right)$.

A \emph{dimension group} is an ordered Abelian group which 
is (order isomorphic to) the inductive limit of a sequence of 
ordered Abelian groups 
\[ \mathbb{Z}^{n_{1}} 
\overset{\alpha_{1}} \longrightarrow 
\mathbb{Z}^{n_{2}} 
\overset{\alpha_{2}} \longrightarrow 
\mathbb{Z}^{n_{3}} 
\overset{\alpha_{3}} \longrightarrow 
\cdots \] 
for some positive integers $n_{j}$ and 
some positive group homomorphisms $\alpha_{j}$.

\medskip

A \emph{supernatural number} is a sequence 
$N = \left\{ n_{j} \right\}^{\infty}_{j = 1}$ 
where each $n_{j}$ belongs to $\{0, 1, 2, \ldots, \infty \}$. 
More suggestively, if $ \left\{ p_1 ,p_2, \ldots \right\}$ is 
the set of all prime numbers listed in increasing order, then 
we may view $N$ as a formal infinite prime factorization
$\prod^{\infty}_{j = 1} p_{j}^{n_{j}}$. 
Then each natural is a supernatural number whose 
sequence is eventually zero. The \emph{product} of 
two supernatural numbers $N = \{ n_{j} \}_{j = 1}^{\infty}$ 
and $M = \{ m_{j} \}_{j = 1}^{\infty}$ is defined to be 
$NM = \{ n_{j} + m_{j} \}_{j = 1}^{\infty}$. 
Also, we write $M \vert N$ if $m_{j} \leq n_{j}$ for all 
$j \geq 1$.

\medskip

The subgroup $\mathit{Q}(N)$ of the additive group 
$\mathbb{Q}$ associated to a supernatural number 
$N = \{ n_{j} \}_{j = 1}^{\infty}$ consists of all fractions 
$x / y$ where $x$ is any integer and 
$y = \prod_{j = 1}^{\infty} p_{j}^{m_{j}}$ 
for some nonnegative integers $m_{j} \leq n_{j}$ where 
$m_{j} > 0$ for only finitely many $j$. Note that the group 
$\mathit{Q}(N)$ is generated by 
\[ \left\{ \frac{1}{p_{1}^{n_{1}}}, \frac{1}{p_{2}^{n_{2}}}, 
\ldots, \frac{1}{p_{k}^{n_{k}}}, \ldots \right\} \] 
in which if $n_{j} = \infty$ for some $j$, then by 
$1/p_{j}^{n_{j}}$ we mean the sequence 
$1/p_{j}, 1/p_{j}^{2}, \ldots$.

We recall the supernatural number $N$ associated to a UHF~algebra $A$.

\begin{definition}[\cite{rll00}]\label{def_supnum}
Let $A$ be a UHF~algebra, that is, a C*-algebra isomorphic to 
the inductive limit of a sequence
\[ \mathrm{M}_{k_{1}} 
\overset{\varphi_{1}} \longrightarrow 
\mathrm{M}_{k_{2}} 
\overset{\varphi_{2}} \longrightarrow 
\mathrm{M}_{k_{3}} 
\overset{\varphi_{3}} \longrightarrow 
\cdots \] 
where the connecting maps $\varphi_{i}$ are unital and where 
$\{ k_{i} \}$ is a sequence of positive integers satisfying 
$k_{i} \vert k_{i+1}$ for all $i \geq 1$. We write 
\[ k_{i} = \prod_{j = 1}^{\infty} p_{j}^{n_{i, j}}, 
\hspace{0.5 cm} n_{i, j} \in \mathbb{Z}^{+}, \] 
and let $N$ be the supernatural number 
$\{ n_{j} \}^{\infty}_{j = 1}$ where 
$n_{j} = \sup \left\{ n_{i, j} : i \in \mathbb{N} \right\}$.
Conversely, if $N=\{ n_{j} \}^{\infty}_{j = 1}$ is 
a supernatural number and we define
\[ \ell_j =\prod\limits_{i=1}^{j} p_{i}^{\min\{ j, n_i \}} \]
for $j \geq 1$, then $\ell_j | \ell_{j+1}$. We denote by 
$\mathrm{M}_N$ the UHF~algebra which is the direct limit of 
$\mathrm{M}_{\ell_j}$'s with the diagonal homomorphisms 
$\varphi_j : \mathrm{M}_{\ell_j} \to \mathrm{M}_{\ell_{j +1}}$ 
as connecting maps. Then $N$ is the supernatural number 
associated to $\mathrm{M}_N$.
\end{definition} 

In the following lemma we gather known facts about 
UHF~algebras needed in the sequel.

\begin{lemma}\label{lem_uhf}
Let $A$ and $B$ be two UHF~algebras with 
supernatural numbers $N = \{ n_{j} \}_{j = 1}^{\infty}$ and 
$M= \{ m_{j} \}_{j = 1}^{\infty}$, respectively.
\begin{enumerate}
\item\label{lem_uhf_it1}  
The following statements are equivalent:
\begin{enumerate}
\item\label{lem_uhf_it1_1} 
There is a unital $*$-homomorphism from $A$ into $B$;

\item\label{lem_uhf_it1_2}
$N \vert M$;

\item\label{lem_uhf_it1_3} 
$\mathit{Q}(N) \subseteq \mathit{Q}(M)$; 

\item\label{lem_uhf_it1_4}
There is a unital (injective) group homomorphism 
from $\mathit{Q}(N)$ into $\mathit{Q}(M)$.
\end{enumerate}

\item\label{lem_uhf_it2}
$A$ is isomorphic to $B$ if and only if there are 
unital $*$-homomorphisms $A \to B$ and $B \to A$.
\end{enumerate}
\end{lemma}

\begin{proof}
We prove \eqref{lem_uhf_it1}.
The equivalence of \eqref{lem_uhf_it1_1}, 
\eqref{lem_uhf_it1_2}, and \eqref{lem_uhf_it1_3} is known 
(see, for instance, \cite[Exercise~7.11]{rll00}). We show that 
\eqref{lem_uhf_it1_3} and \eqref{lem_uhf_it1_4} are 
equivalent.

First let us point out a fact: every unital group homomorphism 
from $\mathit{Q}(N)$ into $\mathit{Q}(M)$ is injective. 
For this, let $\theta : \mathit{Q}(N) \rightarrow \mathit{Q}(M)$ 
be such a homomorphism and let $m/n$ be in $\mathit{Q}(N)$ 
with $\theta(m/n) = 0.$ If $m \neq 0$ then $\theta(1/n) = 0$ 
and hence $0 = n \theta(1/n) = \theta(1) = 1$ that is 
impossible. Thus $m = 0$ and so $m/n = 0$.

Now let $\theta : \mathit{Q}(N) \rightarrow \mathit{Q}(M)$ be 
an injective group homomorphism with $\theta(1) = 1$. 
For every $k, j \in \mathbb{N}$ with 
\( k \leq n_{j} \), we get \( p_{j}^{k} \theta(1/p_{j}^{k}) = \theta(1) = 1 \)
and hence $1/p_{j}^{k}$ belongs to $\mathit{Q}(M)$. Thus
$\mathit{Q}(N) \subseteq \mathit{Q}(M)$. For the converse, 
consider the canonical injection from $\mathit{Q}(N)$ into 
$\mathit{Q}(M)$.

Part~\eqref{lem_uhf_it2} follows from 
Part~\eqref{lem_uhf_it1} and \cite[Proposition~7.4.5]{rll00}.  
\end{proof}

\begin{remark}\label{rem_mu_uniq}
Let $A$ be a unital C*-algebra. Then by 
Lemma~\ref{lem_uhf}\eqref{lem_uhf_it2},
a maximal UHF~subalgebra of $A$ in the sense of
Definition~\ref{def_mu} is unique up to isomorphism (if exists).
Also, by Lemma~\ref{lem_uhf}\eqref{lem_uhf_it1},
a unital UHF~subalgebra $B \cong \mathrm{M}_N$ of $A$ is 
a maximal UHF~subalgebra if $m \vert N$ for any other
unital UHF~subalgebra $D \cong \mathrm{M}_{m}$ of $A$.
\end{remark}

\subsection{Property~(D)}\label{subsec_propd}
In this subsection we introduce Property~(D).
 
\begin{definition}\label{def_vert}
Let $(G, G^{+}, u)$ be 
an ordered Abelian group with distinguished order unit $u$.
\begin{enumerate}
\item\label{def_vert_it1}
If $n$ is a natural number, we write $n \vert u$, if there exists 
$x$ in $G^{+}$ such that $nx = u$.

\item\label{def_vert_it2}
If $N$ is a supernatural number, we write $N \vert u$, if 
$n \vert u$ for all natural numbers $n$ for which $n \vert N$.
\end{enumerate}
\end{definition}  

Note that if $M \vert N$ and $N \vert u$, then $M \vert u$.

\begin{definition}\label{def_propd}
We say that an ordered Abelian group with order unit 
$(G, G^{+}, u)$ has \emph{Property~(D)} if 
every co-prime natural numbers $n$ and $m$ with $n \vert u$ 
and $m \vert u$ satisfy $nm \vert u$. 
\end{definition}

\begin{lemma}\label{lem_weakunp_propd}
Every weakly unperforated ordered Abelian group with 
order unit has Property~(D). In particular, every 
dimension group has Property~(D).
\end{lemma}

\begin{proof}
Let $(G, G^{+}, u)$ be 
a weakly unperforated ordered Abelian group with 
a distinguished order unit. Let $n$ and $m$ be 
co-prime natural numbers such that $n \vert u$ and 
$m \vert u$. So there are $x, y \in G^{+}$ such that 
$mx = ny = u$. Since $\gcd(n, m) = 1$, there are 
$k, l \in \mathbb{Z}$ satisfying $km + l n = 1$. Thus 
$nm (l x + ky) = u$. In particular, $lx + ky >0$ as $G$ is 
weakly unperforated. Therefore, $nm \vert u$.
\end{proof}

Note that the notion of Property~(D) in 
Definition~\ref{def_propd} depends on the order unit. For 
example, consider the positive cone $C = \{ 0, 2, 3, \ldots \}$
for $\mathbb{Z}$. Then $(\mathbb{Z}, C, 2)$ has 
Property~(D), but $(\mathbb{Z}, C, 6)$ does not, since 
$2 \vert 6$ and $3 \vert 6$ but 6 does not divide 6 in this 
ordered Abelian group. Also, note that $(\mathbb{Z}, C, 6)$ is 
not weakly unperforated as $6\cdot 1 > 0$ but $1 \not> 0$.

\medskip

Now we give an equivalent condition to Property~(D). If 
$\Sigma$ is a family of supernatural numbers, by 
the \emph{maximum element} of $\Sigma$ we mean 
the maximum element of the partially ordered set 
$(\Sigma, \precsim)$ where $M \precsim N$ means $M \vert N$.

\begin{theorem}\label{thm_max_propd}
An ordered Abelian group with order unit $(G, G^{+}, u)$ has 
Property~(D) if and only if the set $\Sigma$ of 
supernatural numbers $N$ with $N \vert u$ has 
the maximum element.
\end{theorem}

\begin{proof}

Suppose that $N = \{ k_{j} \}_{j \in \mathbb{N}}$ is 
the maximum element of $\Sigma$, and let $n$ and $m$ be 
co-prime natural numbers such that $n \vert u$ and $m \vert u$.
We can consider $n = \{ n_{j} \}^{\infty}_{j = 1}$ and 
$m = \{ m_{j} \}^{\infty}_{j = 1}$ as supernatural where 
these sequences are eventually zero. Since $n, m \in \Sigma$, 
we see that $n \vert N$ and $m \vert N$. Since 
$\gcd (n, m) = 1$ and $nm \vert N$, $nm \vert u$. Thus 
$(G, G^{+}, u)$ has Property~(D).

\medskip

For the converse, set 
$k_{j} := \sup \{ k \geq 0 : p_{j}^{k} \vert u \}$ and define 
the supernatural number 
$N := \{ k_{j} \}_{j \in \mathbb{N}}$. We show that $N$ 
is in $\Sigma$ and is its maximum. Let 
a natural number $n = p_{1}^{n_{1}} \cdots p_{t}^{n_{t}}$ 
satisfy $n \vert N$. Since $n_{j} \leq k_{j}$, 
$p_{j}^{n_{j}} \vert u$, for all $1 \leq j \leq t$, and hence 
$n \vert u$ as $G$ has Property~(D). By 
Definition~\ref{def_vert}\eqref{def_vert_it2}, 
$N \vert u$ and so $N\in\Sigma$. Finally, let 
$M = \{ l_{j} \}_{j \in \mathbb{N}}$ be in $\Sigma$. For any 
$j$ since $p_{j}^{l_{j}} \vert M$ and $M \vert u$, we get 
$l_{j} \leq k_{j}$ and hence $M \vert N$. Thus $N$ is 
the maximum element of $\Sigma$.
\end{proof}

We denote by $N(G,u)$ the maximum supernatural number dividing 
$u$ defined in the preceding proof.

\subsection{$K_{0}$-lifting property for UHF~algebras}\label{subsec_k0lift}

\begin{definition}\label{def_k0lift}
We say that a unital C*-algebra $A$ has 
\emph{$K_{0}$-lifting property for UHF algebras} if for any 
UHF~algebra $D$, the existence of 
an injective positive order unit preserving group homomorphism 
$K_{0}(D) \to K_{0}(A)$ implies the existence of 
a (necessarily injective) unital $*$-homomorphism $D \to A$. 
\end{definition}

We give a list of C*-algebras having $K_{0}$-lifting property for 
UHF~algebras.

\begin{proposition}\label{prop_list_K0lift}
The following classes of C*-algebras have 
$K_{0}$-lifting property for UHF~algebras:

\begin{enumerate}
\item\label{prop_list_K0lift_it1}
unital AF~algebras,

\item\label{prop_list_K0lift_it2}
unital simple separable C*-algebras with tracial rank zero,

\item\label{prop_list_K0lift_it3}
unital properly infinite C*-algebras.
\end{enumerate}
\end{proposition}

\begin{proof}
Part \eqref{prop_list_K0lift_it1} is known, in fact, if $A$ and 
$D$ are unital AF~algebras and 
$\alpha : K_{0}(D) \to K_{0}(A)$ is 
a positive group homomorphism with 
$\alpha([1_{D}]) = [1_{A}]$, then there is 
a unital $*$-homomorphism $\varphi : D \to A$ such that 
$K_{0}(\varphi) = \alpha$ (see, e.g., 
\cite[Exercise~7.7]{rll00}).
  
Part~\eqref{prop_list_K0lift_it2} follows from \cite[Theorem~6.4]{da04} 
which says that if $D$ and $A$ 
are unital simple separable C*-algebras with tracial rank zero 
such that $D$ is exact and satisfies the UCT, then for 
any $\alpha \in KK(D, A)$ with 
$\alpha_{*}(K_{0}(D)^{+}) \subseteq K_{0}(A)^{+}$ and 
$\alpha_{*}[1_{D}] = [1_{A}]$ there is (up to 
approximately unitarily equivalence) 
a nuclear unital $*$-homomorphism $\varphi: D \to A$ 
such that $\varphi_{*}(x) = \alpha_{*}(x)$ for all 
$x \in \underline{K}(D)$.

Part~\eqref{prop_list_K0lift_it3} follows from \cite[Lemma~7.2]{ro95} 
stating that if $A$ is 
a properly infinite unital C*-algebra and $D$ is 
a unital AF~algebra, then for any group homomorphism
$\alpha : K_{0}(D) \to K_{0}(A)$ with 
$\alpha([1_{D}]) = [1_{A}]$ there is 
a unital $*$-homomorphism $\varphi : D \to A$ such that 
$K_{0}(\varphi) = \alpha$.
\end{proof}

\begin{example}\label{exa_cuntz_jisu_k0lift}
The Cuntz algebras $\mathcal{O}_{n}$ for $2 \leq n \leq \infty$
have $K_{0}$-lifting property for UHF~algebras, by 
Part~\eqref{prop_list_K0lift_it3} of the preceding proposition.

Also, the Jiang-Su algebra $\mathcal{Z}$ has this property, 
however, it is not covered by Proposition~\ref{prop_list_K0lift}.
In fact, let $D$ be a UHF~algebra and 
$\alpha : K_{0}(D) \to K_{0}(\mathcal{Z}) \cong \mathbb{Z}$ 
be an injective positive order unit preserving  homomorphism. 
Consider the natural unital map 
$\iota : \mathbb{C} \to \mathcal{Z}$ and 
the induced isomorphism 
$K_{0}(\iota) : K_{0}(\mathbb{C}) \to K_{0}(\mathcal{Z})$. 
Since 
$K_{0}(\iota)^{-1} \circ \alpha : K_{0}(D) \to K_{0}(\mathbb{C})$ 
is an injective positive order unit preserving homomorphism,
applying Proposition~\ref{prop_list_K0lift}\eqref{prop_list_K0lift_it1}
to the C*-algebra $\mathbb{C}$, we get 
a unital $*$-homomorphism $\varphi : D \to \mathbb{C}$, and 
so $\iota \circ \varphi : D \to \mathcal{Z}$
is the desired homomorphism. It follows also that 
$D \cong \mathbb{C}$.
\end{example}

The following observation enables us to find more C*-algebras 
having $K_{0}$-lifting property for UHF~algebras.

\begin{proposition}\label{prop_k0lift_shift}
Let $A$ and $B$ be unital C*-algebras. Suppose that there are 
a unital $*$-homomorphism $\varphi : A \to B$ and 
an injective positive order unit preserving homomorphism
$\beta : K_{0}(B) \to K_{0}(A)$. If $A$ has 
$K_{0}$-lifting property for UHF~algebras, then so does $B$.
\end{proposition}

\begin{proof}
Let $A$ have $K_{0}$-lifting property for UHF~algebras and let  
$\alpha : K_{0}(D) \to K_{0}(B)$ be
an injective positive order unit preserving group homomorphism
for some UHF~algebra $D$. Consider 
the injective positive order unit preserving homomorphism 
$\beta \circ \alpha : K_{0}(D) \to K_{0}(A)$. Then we get 
a unital $*$-homomorphism $\eta: D \to A$, and so 
$\varphi \circ \eta : D \to B$ is 
the desired unital $*$-homomorphism.
\end{proof}
 
\begin{corollary}\label{cor_hequi_k0lift}
Let $A$ and $B$ be unital C*-algebras with $A \sim_h B$. Then 
$A$ has $K_{0}$-lifting ptoperty for UHF~algebras if and only if
so does $B$.
\end{corollary}

\begin{proof}
Let 
$B \overset{\psi} \longrightarrow A \overset{\varphi} \longrightarrow B$
be a homotopy between $A$ and $B$. Then $\varphi$ and 
$\psi$ are unital. For this, note that $\psi \circ \varphi(1_{A})$ 
and $1_{A}$ are the homotopy equivalent projections, and 
hence they are unitarily equivalent. Thus 
$\psi \circ \varphi(1_{A}) = 1_{A}$, and
$\varphi \circ \psi(1_{B}) = 1_{B}$. Since $\varphi(1_{A})$ and
$\psi(1_{B})$ are projections in $B$ and $A$, respectively, 
we get $\varphi(1_{A}) \leq 1_{B}$ and 
$\psi(1_{B}) \leq 1_{A}$, and hence 
$1_{B} = \varphi \circ \psi(1_{B}) \leq \varphi(1_{A}) \leq 1_{B}$. 
Thus $\varphi(1_{A}) = 1_{B}$. Similarly, 
$\psi(1_{B}) = 1_{A}$.

By \cite[Proposition~3.2.6]{rll00},
$K_0 (\psi) :  K_{0}(B) \to K_{0}(A)$ is an isomorphism. Then 
Proposition~\ref{prop_k0lift_shift} implies the statement. 
\end{proof}

\medskip

As an example, for any contractible compact Hausdorff space 
$X$, the C*-algebra $C(X)$ has 
$K_{0}$-lifting ptoperty for UHF~algebras as 
$C(X) \sim_h \mathbb{C}$.
 
\medskip
 
As another application of Proposition~\ref{prop_k0lift_shift},
if $A$ and $B$ are unital C*-algebras and $A \oplus B$ has 
$K_{0}$-lifting property for UHF~algebras, then so do have 
both $A$ and $B$. Also, we have the following result.

\begin{corollary}\label{cor_exaseq_k0lift}
Let there exist a split exact sequence 
\[ 0 \longrightarrow 
I \longrightarrow 
A \overset{\varphi}{\underset{\psi} \rightleftarrows} 
B \longrightarrow 
0 \] 
where $I$ is a C*-algebra, $A$ and  $B$ are 
unital C*-algebras, and $\varphi, \psi$ are 
unital $*$-homomorphisms. If $A$ has 
$K_{0}$-lifting property for UHF~algebras then so does $B$.
\end{corollary}

\subsection{Proof of Theorem~\ref{thm_main}} \label{subsec_proofb}
Let $A$ be C*-algebra with 
$K_{0}$-lifting property for UHF~algebras such that 
$(K_{0}(A), K_{0}(A)^{+})$ is unperforated.

\emph{Existence}: 
Let $\Sigma$ be the set of all supernatural numbers $m$ 
such that $m \vert [1]_{0}$. Then by 
Lemma~\ref{lem_weakunp_propd} and 
Theorem~\ref{thm_max_propd}, $\Sigma$ has 
the maximum element $N = \{ n_{j} \}_{j=1}^{\infty}$ where
\( n_{j} = \sup \{ k \geq 0 : p_{j}^{k} | [1]_{0} \} \) for all 
\( j \in \mathbb{N} \). Consider the UHF~algebra 
$\mathrm{M}_{N}$ and note that 
$(K_{0}(\mathrm{M}_{N}), [1]_0) \cong (\mathit{Q}(N), 1)$ 
as ordered groups with distinguished order unit. We show that 
$K_{0}(\mathrm{M}_{N})$  embeds into $K_{0}(A)$. Let 
$\mathit{Q}(N) = \cup_{j=1}^{\infty} \ell_{j}^{-1} \mathbb{Z}$ 
where $\ell_{j}$ is as in Definition~\ref{def_supnum}. For any
$j \in \mathbb{N}$, since $\ell_{j} \vert N$ and 
$N \vert [1]_{0}$, there is $x_{j} \in K_{0}(A)^{+}$ such that 
$\ell_{j} x_{j} = [1]_{0}$. Now we define 
a positive order preserving group homomorphism 
\begin{align*}
\alpha : 
&\mathit{Q}(N) \to K_{0}(A) \\
&k/\ell_{j} \mapsto k x_{j} 
\end{align*}
where $j \in \mathbb{N}$ and $k \in \mathbb{Z}$. First 
we show that $\alpha$ is well defined. For 
$j, j^{\prime} \in \mathbb{N}$ with $j < j^{\prime}$ and 
$k, k^{\prime} \in \mathbb{Z}$, let 
$k/\ell_{j} = k^{\prime}/\ell_{j^{\prime}}$. Since
$\ell_{j} \left(x_{j} - 
(\ell_{j^{\prime}}/\ell_{j}) x_{j^{\prime}}\right) = 0$ 
and
$K_0 (A)$ is torsion-free (as it is unperforated),  
$x_{j} - (\ell_{j^{\prime}}/\ell_{j}) x_{j^{\prime}} = 0$ and
$x_{j} = (\ell_{j^{\prime}}/\ell_{j}) x_{j^{\prime}}$. Thus
$k^{\prime} x_{j^{\prime}} = 
k (\ell_{j^{\prime}}/\ell_{j}) x_{j^{\prime}} =k x_{j}$, 
as desired. Since $K_0 (A)$ is torsion-free, $\alpha$ is injective.

By assumption, $A$ has 
$K_0$-lifting property for UHF~algebras, and so there is 
a unital $*$-homomorphism 
$\varphi : \mathrm{M}_{N} \to A$. Set 
$MU(A) := \varphi(\mathrm{M}_{N}).$

\medskip

\emph{Maximality}: Let $D \cong \mathrm{M}_{m}$ be 
a unital UHF~subalgebra of $A$ with 
$m = \{ m_{j} \}_{j=1}^{\infty}$. Consider 
the homomorphism $K_{0}(\iota): K_{0}(D) \to K_{0}(A)$ 
where $\iota : D \to A$ is the canonical injection. For any 
natural numbers $j$ and $k \leq m_{j}$, 
$p_{j}^{k} \vert [1]_{0}$ in $K_0(D)$ as
$K_{0}(D) \cong Q(m)$. Thus $p_{j}^{k} \vert [1]_{0}$
in $K_{0}(A)$. Hence $m | [1]_0$ and so $m \in \Sigma$. Thus
$m \vert N$. Therefore, $D$ embeds into $MU(A)$ by 
Remark~\ref{rem_mu_uniq}.
\qed

\begin{corollary}\label{cor_equ_cond}
Let $A$ be a unital C*-algebra such that $K_{0}(A)$ is unperforated. 
Let $N=N\left(K_{0}(A), [1]_0\right)$ be 
the maximum supernatural number as in 
Theorem~\ref{thm_max_propd}. Then the following are 
equivalent:
\begin{enumerate}

\item
$A$ has a maximal UHF~subalgebra 
$MU(A) \cong \mathrm{M}_{N}$,
 
\item
$\mathrm{M}_{N}$ embeds unitaly in $A$,
 
\item
$A$ has the $K_{0}$-lifting property for UHF~algebras.
\end{enumerate}
\end{corollary}

\begin{proof}
First note that the maximum supernatural number 
$N(K_{0}(A), [1]_0)$ as in Theorem~\ref{thm_max_propd}
exists since $K_{0}(A)$ is unperforated and 
Lemma~\ref{lem_weakunp_propd} can be applied.

Now for $(1) \Rightarrow (3)$, let $D$ be a UHF~algebra and 
$\alpha : K_{0}(D) \to K_{0}(A)$ be 
a positive injective order unit preserving homomorphism. 
It follows that $N_{D} \vert N$, and hence by 
Lemma~\ref{lem_uhf}, there is a unital $*$-homomorphism 
$\psi : D \to MU(A)$. Then $\iota_{MU(A)} \circ \psi : D \to A$
is a unital $*$-homomorphism. Thus $A$ has 
the $K_{0}$-lifting property for UHF~algebras. 

For $(2) \Rightarrow (1)$, let 
$\varphi : \mathrm{M}_{N} \to A$ be a unital embedding. Then
by the part ``maximality" of the proof of  
Theorem~\ref{thm_main}, $\varphi(\mathrm{M}_{N})$ is  
a maximal UHF~subalgebra of $A$. 

For $(3) \Rightarrow (2)$, by Theorem~\ref{thm_main}, $A$ 
has a maximal UHF~subalgebra. The part ``existence" of 
the proof of Theorem~\ref{thm_main} implies that 
$MU(A) \cong \mathrm{M}_{N}$.
\end{proof}

\medskip

There are examples of unital C*-algebras $A$ which have 
a maximal UHF~subalgebra but $K_{0}(A)$ is perforated. For 
example, $K_{0}\left(C(\mathbb{T}^{4})\right) \cong \mathbb{Z}^{8}$
is perforated (by \cite[Example~6.7.2(b)]{bl86} and 
\cite{el84}), however, 
$MU\left(C(\mathbb{T}^{4})\right) \cong \mathbb{C}$ (by 
Proposition~\ref{prop_mu} below).

\begin{remark}\label{rem_mu_embed}
In Definition~\ref{def_mu}, if we do not assume that  
$1_{D} = 1_{A}$ then unusual examples arise. For instance, 
consider the AF~algebra $A= \mathcal{K} + \mathbb{C}1$ and
its maximal UHF~subalgebra $MU(A)$ in the sense of 
this new definition. Since for any $m \geq 1$, the matrix algebra
$\mathrm{M}_{m}$ embeds into 
$\mathcal{K} + \mathbb{C}1$ (by a nonunital embedding), it 
also embeds into $MU(A)$. However, every 
unital simple C*-subalgebra $B$ of $A$ is finite dimensional. 
In fact, If $B \cap \mathcal{K} \neq \{0\}$, then 
$B \cap \mathcal{K} = B$ (since 
$B \cap \mathcal{K}  \unlhd B$). As 
$\mathcal{K}$ is liminal, so is $B$, and therefore 
$B$ is finite dimensional. Now let 
$B \cap \mathcal{K} = \{0\}$. If $1_{B} = 1_{A}$ 
then it follows that $B=\mathbb{C}1_{A}$. If 
$1_{B} \neq 1_{A}$, then there is a nonzero projection 
$p \in \mathcal{K}$ such that $1_{B} =1_{A} - p$, 
and it follows that $B = \mathbb{C}(1-p)$. Therefore, such 
a maximal UHF~subalgebra $MU(A)$ of $A$ does not exist.

As another example, if 
$A = \mathrm{M}_{4} \oplus \mathrm{M}_{6}$ and 
its maximal UHF~subalgebra $MU(A)$ in the sense of 
this new definition exists, then $\mathrm{M}_{4}$ and 
$\mathrm{M}_{6}$ embed unitaly into $MU(A)$. Now it follows
from Lemma~\ref{lem_uhf}\eqref{lem_uhf_it1} that 
$\mathrm{M}_{12}$ embeds into $MU(A)$, and hence into 
$A$, which is impossible.
\end{remark}

\section{Permanence properties and examples of  
a maximal UHF~subalgebra}\label{sec_exa}

The following results enable us to find examples of  
a maximal UHF~subalgebra. The proof of the first one is 
a direct application of Definition~\ref{def_mu} and so 
is omitted.

\begin{proposition}\label{prop_mu_homo}
Let $A$ and $B$ be unital C*-algebras. Let there be 
unital $*$-homomorphisms $\varphi : A \to B$ and 
$\psi : B \to A$. Then $A$ has a maximal UHF~algebra 
if and only if so does $B$. In this case, $MU(A) \cong MU(B)$. 
\end{proposition}

The following result is about split exact sequences: 

\begin{corollary}\label{cor_mu_exactseq}
With the assumptions of Corollary~\ref{cor_exaseq_k0lift}, if 
moreover $A$ has 
the $K_{0}$-lifting property for UHF~algebras and $K_{0}(A)$ 
is unperforated, then $B$ has a maximal UHF~algebra, and 
$MU(A) \cong MU(B)$.  
\end{corollary}

\begin{corollary}\label{cor_mu_hequi}
Let $A$ and $B$ be unital C*-algebras with $A \sim_{h} B$. If 
$A$ has a maximal UHF~subalgebra then so does $B$ and 
$MU(A) \cong MU(B)$.
\end{corollary}

\begin{proof}
The proof of Corollary~\ref{cor_hequi_k0lift} provides 
unital $*$-homomorphisms $\varphi : A \to B$ and 
$\psi : B \to A$. Then Proposition~\ref{prop_mu_homo} 
can be applied.
\end{proof}

As an example, if $X$ is a compact Hausdorff contractable space 
and $A$ is a unital C*-algebra which has 
a maximal UHF~subalgebra, then so does $A \otimes C(X)$ and 
$MU\left(A \otimes C(X)\right) \cong MU(A)$. This follows from 
the fact that $C(X) \sim_{h} \mathbb{C}$ and so 
$A \otimes C(X) \sim_{h} A$.

\begin{lemma}\label{lem_max_ind_lim}
Let 
\[ A_{1} 
\overset{\varphi_{1}} 
\longrightarrow 
A_{2}
\overset{\varphi_{2}} 
\longrightarrow 
A_{3}
\overset{\varphi_{3}} 
\longrightarrow 
\cdots
\longrightarrow 
A, \] 
be an inductive limit of unital C*-algebras $A_{n}$ such that
every \( A_{n} \) has a maximal UHF~subalgebra and 
unital connecting maps $\varphi_{n}$, and
$MU(A_{n}) \cong \mathrm{M}_{N_{n}}$ for 
a supernatural number \( N_{n} \). If $K_{0}(\varphi_{n})$ 
is injective for all $n$, then $A$ has a maximal UHF~subalgebra
$MU(A)$ that is isomorphic to 
$\mathrm{M}_{\sup\limits_{n \in \mathbb{N}} N_{n}}$. 
In other words, 
$MU(\underrightarrow{\lim} A_{n}) \cong  
\underrightarrow{\lim} MU(A_{n})$.
\end{lemma}

\begin{proof}
Let $\Sigma$ be the set of all supernatural numbers $m$ 
such that $m \vert [1_{A}]_{0}$. using 
\cite[Theorem~6.3.2(ii)]{rll00} and since for all $n$, 
$K_{0}(\varphi_{n})$ is injective, it follows that $\Sigma$ has 
the maximum element 
$N = \sup_{n \in \mathbb{N}} N_{n}$. We will show that
$\mathrm{M}_{N}$ embeds into $A$. 

Let for a UHF~algebra $\mathrm{M}_{K}$, there is 
a unital $*$-homomorphism from $\mathrm{M}_{K}$ into $A$.
Hence $K \vert [1_{A}]_{0}$ and using
\cite[Theorem~6.3.2(ii) and (iii)]{rll00}, we see that $K \leq N$. 
According to Lemma~\ref{lem_uhf}, there is a unital embedding 
from $\mathrm{M}_{K}$ into $\mathrm{M}_{N}$. Thus by 
Definition~\ref{def_mu}, $A$ has a maximal UHF~subalgebra 
$MU(A)$ and $MU(A) \cong \mathrm{M}_{N}$. Note that
for any \( n \in \mathbb{N} \) there is an embedding 
\( \psi_{n} : \mathrm{M}_{N_{n}} \to \mathrm{M}_{N_{n+1}} \)
and \( \mathrm{M}_{N} \) is isomorphic to 
the resulting inductive sequence:
\[ \mathrm{M}_{N_{1}} 
\overset{\psi_{1}} 
\longrightarrow 
\mathrm{M}_{N_{2}} 
\overset{\psi_{2}} 
\longrightarrow 
\mathrm{M}_{N_{3}} 
\overset{\psi_{3}} 
\longrightarrow 
\cdots
\longrightarrow 
\mathrm{M}_{N} . \] 
If fact, \( \varphi_{n}(MU(A_{n})) \) is isomorphic to 
\( MU(A_{n}) \) and hence it is a unital UHF~subalgebra of 
\( A_{n+1} \). By Definition~\ref{def_mu}, there is an embedding 
\( \theta_{N} : \varphi_{n}(MU(A_{n})) \to MU(A_{n+1}) \).
Then \( \psi_{n} : MU(A_{n}) \to MU(A_{n+1}) \) denoted by 
\( \psi_{n}(x) = \theta_{n} \circ \varphi_{n}(x) \), is 
the desired connecting map. (Note that the inductive limit of 
\( \mathrm{M}_{N_{n}} \)'s is (up to isomorphism) independent of
\( \psi_{n} \)'s and is isomorphic to \( \mathrm{M}_{N} \).)
\end{proof}

\begin{corollary}\label{cor_mu_oplus}
Let $A$ and $B$ be unital C*-algebras. If $A \oplus B$ has 
the $K_{0}$-lifting property for UHF~algebras and 
$K_{0}(A \oplus B)$ is unperforated, then $A$ and $B$ have  
maximal UHF~subalgebras.
\end{corollary}

\begin{proof}
First note that, for unital C*-algebras $A$ and $B$, 
$K_{0}(A \oplus B)$ is unperforated if and only if so are 
both $K_{0}(A)$ and $K_{0}(B)$. Now, let $A$ and $B$ be 
as in the statement. Then $K_{0}(A)$ and $K_{0}(B)$ are 
unperforated. On the other hand, by the remark preceding
Corollary~\ref{cor_exaseq_k0lift}, both $K_{0}(A)$ and 
$K_{0}(B)$ have the $K_{0}$-lifting property for UHF~algebras.
Finally, applying Theorem~\ref{thm_main} we get the result.
\end{proof}

In the following proposition we give a list of examples of 
maximal UHF~subalgebras of certain C*-algebras.

\begin{proposition}\label{prop_mu}
\hfill
\begin{enumerate}
\item\label{prop_mu_it1}
A maximal UHF~subalgebra of any 
finite dimensional C*-algebra 
$\mathrm{M}_{k_{1}} \oplus \cdots \oplus \mathrm{M}_{k_{r}}$ 
is isomorphic to $\mathrm{M}_{\gcd (k_{1}, \ldots, k_{r})}$;

\item\label{prop_mu_it2}
For any UHF~algebra $A$ and any 
compact Hausdorff contractible space $X$, 
$MU\left(A \otimes C(X)\right) \cong A$;

\item\label{prop_mu_it3}
a maximal UHF~sabalgebra of the following C*-algebras is 
isomorphic to $\mathbb{C}$:

\begin{enumerate}
\item\label{prop_mu_it31}
every unital projectionless C*-algebra;

\item\label{prop_mu_it32}
every unital C*-algebra with a projectionless quotient;

\item\label{prop_mu_it33}
every unital C*-algebra which is not divisible;

\item\label{prop_mu_it34}
the unitization algebra $\tilde{A}$ of any 
nonunital C*-algebra $A$;

\item\label{prop_mu_it35}
every unital C*-algebra having a character; 

\item\label{prop_mu_it36}
the Toeplitz algebra $\mathcal{T}$;

\item\label{prop_mu_it37}
every unital Abelian 
C*-algebra;

\item\label{prop_mu_it38}
the unital universal C*-algebra generated by two projections;

\item\label{prop_mu_it39}
the Cuntz algebra $\mathcal{O}_{\infty}$;
\end{enumerate}

\item\label{prop_mu_it4} 
a maximal UHF~sabalgebra of the following C*-algebras is 
isomorphic to the universal UHF~algebra $\mathcal{Q}$:

\begin{enumerate}
\item\label{prop_mu_it41}
$B(H)$ for any infinite dimensional Hilbert space $H$;

\item\label{prop_mu_it42}
the Calkin algebra $\mathcal{Q}(H)$ for any 
infinite dimensional Hilbert space $H$;

\item\label{prop_mu_it43}
$\mathcal{M}(A)$ and $\mathcal{M}(A)/A$ for any 
nonzero stable C*-algebra $A$; 

\item\label{prop_mu_it44}
the Cuntz algebra $\mathcal{O}_{2}$;

\item\label{prop_mu_it45}
every unital  C*-algebra generated by two isometries satisfying 
the Cuntz relation.
\end{enumerate}
\end{enumerate}
\end{proposition}

\begin{proof}
\hfill
\begin{enumerate}
\item
Let 
$A = \mathrm{M}_{k_{1}} \oplus \cdots \oplus \mathrm{M}_{k_{r}}$ 
and $k = \gcd(k_{1}, \ldots, k_{r})$. Since $k \vert k_{j}$ for 
all $1 \leq j \leq r$,  there is a unital $*$-homomorphism 
$\psi_{j} : \mathrm{M}_{k} \to \mathrm{M}_{k_{j}}$. We 
consider the $*$-homomorphism 
$\psi : \mathrm{M}_{k} \to A$ by setting 
$\psi(x) = \left(\psi_{1}(x), \ldots, \psi_{r}(x)\right)$, 
$x \in \mathrm{M}_{k}$, and so we see that 
$\mathrm{M}_{k}$ embeds into $A$. Now let $D$ be 
a UHF~subalgebra of $A$. Thus there is $l \geq 1$ and 
an isomorphism $\theta : \mathrm{M}_l \to D$.
Consider the unital $*$-homomorphism 
$\pi_{j} \circ \iota_{D} \circ \theta : \mathrm{M}_{l} 
\to \mathrm{M}_{k_{j}}$ 
where $\pi_{j}$ is the  projection map from $A$ onto 
$\mathrm{M}_{k_{j}}$. Hence $l \vert k_{j}$ for all 
$1 \leq j \leq r$, and so $l \vert k$. Hence, 
$\mathrm{M}_{l}$ embeds into $\mathrm{M}_{k}$, and 
so does $D$ into $\psi(\mathrm{M}_{k})$. Therefore, 
$MU(A)=\psi(\mathrm{M}_{k})$.

\item
Since $A \otimes C(X) \sim_{h} A$, we see that
$MU\left(A \otimes C(X)\right) \cong MU(A)=A$ (see 
the example following Corollary~\ref{cor_mu_hequi}).

\item
\eqref{prop_mu_it31}, \eqref{prop_mu_it32}, and 
\eqref{prop_mu_it33} are clear. 
Proposition~\ref{prop_mu_homo} implies 
\eqref{prop_mu_it34}, \eqref{prop_mu_it35}, 
\eqref{prop_mu_it36}, \eqref{prop_mu_it37}, and 
\eqref{prop_mu_it38}. See Example~\ref{exa_mu_part2} 
below for \eqref{prop_mu_it39}.

\item
Part~\eqref{prop_mu_it41} follows from the fact that 
every separable unital C*-algebra is embedded unitally into 
$B(H)$. Consider the quotient map from $B(H)$ onto 
$\mathcal{Q}(H)$. Then, using \eqref{prop_mu_it41}, we get 
\eqref{prop_mu_it42}. Part~\eqref{prop_mu_it43} follows from
\eqref{prop_mu_it41} and the fact that $B(\ell^{2})$ 
embeds into both $\mathcal{M}(A)$ and $\mathcal{M}(A)/A$ 
\cite[Paragraph~5.1.9]{li01}. Since every 
separable exact unital C*-algebra is embedded unitally into 
$\mathcal{O}_{2}$ (\cite[Theorem~6.3.11]{rs02}), we get 
\eqref{prop_mu_it44}. Part~\eqref{prop_mu_it45} follows from
\eqref{prop_mu_it44} and the universal property of 
$\mathcal{O}_{2}$ \cite{rs02}.
\qedhere
\end{enumerate}
\end{proof}

According to the following result, 
every simple infinite unital C*-algebra $A$ contains 
a unital subalgebra $B$ such that a maximal UHF~subalgebra of 
a quotient of $B$ is isomorphic to $\mathcal{Q}$. 

\begin{corollary}\label{prop_iso_cuntzn}
For every simple infinite unital C*-algebra $A$ there is 
a unital C*-subalgebra $B$ of $A$ and a closed ideal $J$ of $B$ 
such that $B/J$ has a maximal UHF~subalgebra isomorphic to 
$\mathcal{Q}$.
\end{corollary}

\begin{proof}
By \cite[Paragraph~3.2]{cu77a} and 
\cite[Paragraph~2.2]{cu77b}, every 
simple infinite unital C*-algebra contains isometries 
$V_{1}, V_{2}$ satisfying 
$V_{1} V_{1}^{*} + V_{2} V_{2}^{*} \leq 1$. By 
\cite[Paragraph~3.1]{cu77a}, there is a closed ideal $J$ of 
$C^{*}(V_{1}, V_{2})$ such that $J \cong \mathcal{K}$ and 
the quotient $C^{*}(V_{1}, V_{2}) / J$ is isomorphic to 
$\mathcal{O}_{2}$. Therefore, by 
Proposition~\ref{prop_mu}\eqref{prop_mu_it44}, 
$MU(C^{*}(V_{1}, V_{2}) / J) \cong \mathcal{Q}$.
\end{proof}

Recall that a Kirchberg algebra is 
a purely infinite, simple, nuclear, separable C*-algebra 
\cite[Definition~4.3.1]{rs02}.

\begin{proposition}\label{prop_o2_embed} 
Let $A$ be a C*-algebra with a properly infinite, full projection
$p$ satisfying $[p]_{0} = 0$ in $K_{0}(A)$. Then $pAp$ has  
a maximal UHF~subalgebra isomorphic to $\mathcal{Q}$. 
This is the case, in particular, when $A$ is 
a unital Kirchberg algebra.
\end{proposition}

\begin{proof}
By assumptions and \cite[Proposition 4.2.3(ii)]{rs02}, there is
a unital $*$-homomorphism $\mathcal{O}_{2} \to pAp$.
Now since $MU(\mathcal{O}_{2}) \cong \mathcal{Q}$ 
(Example~\ref{prop_mu}\eqref{prop_mu_it44}), we get 
$MU(pAp) \cong \mathcal{Q}$. Also it is proved in 
\cite[Theorem~4.1]{cu81} that if $A$ contains 
a properly infinite, full projection, then
\[ K_{0}(A) = \left\{ [p]_{0} : p \ \text{is a properly infinite, 
full projection in} \ A \right\}. \]
Thus every unital Kirchberg algebra $A$ has 
a properly infinite, full projection $p$ such that $[p]_{0} = 0$ 
and therefore $MU(pAp) \cong \mathcal{Q}$. 
\end{proof}

In particular, for every Cuntz algebra $\mathcal{O}_{n}$ for 
$3 \leq n \leq \infty$, there is a corner whose 
maximal UHF~subalgebra is isomorphic to $\mathcal{Q}$.

\medskip

Not all unital C*-algebras have a maximal UHF~subalgebra. 
First we need the following lemma essentially contained in 
\cite{rv98}.

\begin{lemma}\label{lem_no_matrix_embed}
Let $k, l \geq 2$ be natural numbers with $k$ prime and $l$ 
not divisible by $k$. Then the universal unital free product 
$\mathrm{M}_{k} \ast_{r} \mathrm{M}_{l}$ does not admit 
any unital embedding of $\mathrm{M}_{kl}$.
\end{lemma}

\begin{proof}
Let 
$\mathcal{A} = \mathrm{M}_{k} \ast_{r} \mathrm{M}_{l}$ 
and
$\tau : \mathcal{A} \rightarrow \mathrm{M}_{k} \otimes \mathrm{M}_{l}$
be the $*$-homomorphism induced by 
the natural $*$-homomorphisms 
$\mathrm{M}_{k} \rightarrow \mathrm{M}_{k} \otimes \mathrm{M}_{l}$ 
and 
$\mathrm{M}_{l} \rightarrow \mathrm{M}_{k} \otimes \mathrm{M}_{l}$,
using the universal property of $\mathcal{A}$. By  
Proposition~3.5 and Theorem~3.6 of \cite{rv98},
$K_{0}(\tau): (K_{0}(\mathcal{A}),K_{0}(\mathcal{A})^{+})  
\to (\mathbb{Z}, \langle k, l \rangle)$ 
is an isomorphism where
$\langle k, l \rangle = \left\{ nk+ml : n, m \in \mathbb{Z}^{+} \right\}$. 
We show that $K_{0}(\tau)([1_{\mathcal{A}}]_{0}) = kl$. Let 
$e_j$ be the matrix in $\mathrm{M}_{k}$ 
having 1 in $jj$-th entry and 0 elsewhere, for $1 \leq j \leq k$. 
Then by 
\cite[Proposition~3.1.7]{rll00},
\begin{align*}
K_{0}(\tau)([1_{\mathcal{A}}]_{0}) 
= [\tau(1_{\mathcal{A}})]_{0} 
&= [1_{M_{k}} \otimes 1_{M_{l}}]_{0} \\
&= [(\sum_{j=1}^{k} e_{j}) \otimes 1_{M_{l}}]_{0} 
= \sum_{j=1}^{k} [e_{j} \otimes 1_{M_{l}}]_{0} \\
&= \sum_{j=1}^{k} \text{rank}(e_{j} \otimes 1_{M_{l}}) 
= \sum_{j=1}^{k} l = kl.
\end{align*}
(See the proof of \cite[Proposition 3.6]{rv98}.) Hence 
we get $K_{0}(\tau)([1_{\mathcal{A}}]_{0}) = kl.$ Now 
we show that there is no unital embedding 
$\varphi : \mathrm{M}_{kl} \rightarrow \mathcal{A}$. 
Suppose that such a map exists. Consider the following diagram
\[ \xymatrix{
(K_{0}(\mathrm{M}_{kl}), K_{0}(\mathrm{M}_{kl})^{+}, 
[1_{\mathrm{M}_{kl}}]_{0}) \ar[d]_{\cong} \ar[r]^{K_{0}(\varphi)} 
&(K_{0}(\mathcal{A}), K_{0}(\mathcal{A})^{+}, 
[1_{\mathcal{A}}]_{0}) \ar[d]_{\cong}^{K_{0}(\tau)} \\
(\mathbb{Z}, \mathbb{Z}^{+}, kl) \ar@{-->}[r]_{\theta}
&(\mathbb{Z}, \langle k, l \rangle, kl)
} \]
where $\theta$ is 
the positive order preserving group homomorphism such that 
the preceding diagram commutes. Since 
$K_{0}(\tau)([1_{\mathcal{A}}]_{0}) = kl$, we get $\theta(1) =1$. 
But $1 \notin \langle k, l \rangle$ and so 
we get a contradiction. Thus $\mathcal{A}$ does not admit 
any embedding of $\mathrm{M}_{kl}$.
\end{proof}

\begin{example}\label{exa_mu_not}
For co-prime numbers $k, l \geq 2$, the unital C*-algebra 
$\mathrm{M}_{k} \ast_{r} \mathrm{M}_{l}$ does not have  
a maximal UHF~subalgebra. For this, let 
$\mathcal{A} = \mathrm{M}_{k} \ast_{r} \mathrm{M}_{l}$ 
and suppose that $MU(\mathcal{A})$ exists and is isomorphic to
a UHF~algebra $\mathrm{M}_{N}$. Since $\mathrm{M}_{k}$ 
and $\mathrm{M}_{l}$ are embedded in 
$MU(\mathcal{A})$, by Lemma~\ref{lem_uhf}, $k \vert N$ and 
$l \vert N$ and hence $kl \vert N$ since $\gcd(k, l) = 1$. Then 
by Lemma~\ref{lem_uhf}, $\mathrm{M}_{kl}$ 
is embedded unitally into $\mathcal{A}$, contradicting 
Lemma~\ref{lem_no_matrix_embed}.
\end{example}

In view of Theorem~\ref{thm_main}, it is natural to search 
for C*-algebras whose maximal UHF~subalgebras are 
isomorphic to a given UHF~algebra $B$. In the rest of 
this section, we do this.

\begin{lemma}\label{lem_irr_iso}
Let $\alpha, \beta$ be distinct irrational numbers and $G, H$ be 
additive subgroups of $\mathbb{Q}$ such that $1 \in G$. If 
$(G + \alpha H, 1) \cong (G + \beta H, 1)$ 
as dimension groups (with order induced from the natural order
on $\mathbb{R}$) with distinguished order unit, then 
$G + \alpha H = G + \beta H$. 
\end{lemma}

\begin{proof}
First, let $\theta : (G + \alpha H, 1) \to (G + \beta H, 1)$ be 
an order isomorphism. Let $x \in G + \alpha H$. Then for any 
nonzero integer $l$, there is an integer $k$ such that $x$   
belongs to the interval $[k/l, (k+1)/l)$. Since $\theta$ is 
order preserving and $\theta(1)=1$, it follows that $\theta(x)$ 
belongs to the same interval. Thus 
$\vert \theta(x)-x \vert \leq 1/l$. Letting $l \to \infty$, we get 
$\theta(x) = x$. Then $G + \alpha H \subseteq G + \beta H$ 
and similarly $G + \beta H \subseteq G + \alpha H$.
\end{proof}

\begin{theorem}\label{thm_mu_af}
For any UHF~algebra $B$ there exists an uncountable family of 
pairwise non-isomorphic simple unital AF~algebras with 
a maximal UHF~subalgebra isomorphic to $B$.
\end{theorem}

\begin{proof}
Let $B \cong \mathrm{M}_{N}$ for a supernatural $N$, and 
consider the simple dimension group 
$\mathit{Q}(N) + \alpha \mathbb{Z}$ where $\alpha$ is 
an arbitrary irrational number. According to  
\cite[Proposition~7.2.8]{rll00}, there is an AF~algebra 
$A(\alpha)$ such that 
$K_{0}\left(A(\alpha)\right) \cong \mathit{Q}(N) + \alpha \mathbb{Z}$ 
as ordered groups. As \( Q(N) + \alpha \mathbb{Z} \) is 
a simple dimension group, the AF algebra \( A(\alpha) \) is simple.
By Proposition~\ref{prop_list_K0lift}\eqref{prop_list_K0lift_it1},   
$A(\alpha)$ has $K_{0}$-lifting property for UHF~algebras. 
Also, $K_{0}(B)\cong Q(N)$ embeds into 
$K_{0}\left(A(\alpha)\right) \cong \mathit{Q}(N) + \alpha \mathbb{Z}$. 
Hence, $B$ embeds unitaly into $A(\alpha)$. By 
Theorem~\ref{thm_main}, a maximal UHF~subalgebra 
$MU(A(\alpha))$ of $A(\alpha)$ exists and is isomorphic to 
$\mathrm{M}_{K}$ for some supernatural number $K$. Using 
Lemma~\ref{lem_uhf}, we see that $N \vert K$. 
The injection map from $MU(A(\alpha))$ into $A(\alpha)$ 
induces a positive homomorphism 
$\theta: \mathit{Q}(K)\to \mathit{Q}(N) + \alpha \mathbb{Z}$. 
Since $\theta(1)=1$, it follows that the range of this map
is contained in $\mathit{Q}(N)$. Hence, by 
Lemma~\ref{lem_uhf}, $K \vert N$ and so $K=N$. Thus 
$MU\left(A(\alpha)\right) \cong B$.

Elementary facts in Linear Algebra imply that there is 
an uncountable set $I$ of irrational numbers such that for 
any distinct $\alpha, \beta \in I$, the set $\left\{ 1, \alpha, \beta \right\}$ 
is $\mathbb{Q}$-linearly independent. Now if 
$\alpha, \beta \in I$ are distinct, then by
Lemma~\ref{lem_irr_iso},
$Q(N) + \alpha \mathbb{Z} \neq Q(N) + \beta \mathbb{Z}$. 
Thus the AF~algebras $A(\alpha)$ and $A(\beta)$ 
are not isomorphic by \cite[Corollary~7.2.11]{mu90}. 
Therefore, $\{ A(\alpha) : \alpha \in I \}$ is the desired family.
\end{proof}

\begin{theorem}\label{thm_mu_taf}
For any UHF~algebra $B$ there exists 
a simple unital separable tracial rank zero algebra $A$ that 
is not an AF~algebra and $MU(A) \cong B$.
\end{theorem}

\begin{proof}
Let $B$ be a UHF~algebra. According to 
\cite[Theorem~A.6]{da04}, there is 
a simple unital separable tracial rank zero C*-algebra $A$  
such that $K_{0}(A) \cong K_{0}(B)$ and 
$K_{1}(A) \cong \mathbb{Z}$. In particular, $A$ is not 
an AF~algebra. By 
Proposition~\ref{prop_list_K0lift}\eqref{prop_list_K0lift_it2}, 
$A$ has $K_{0}$-lifting property for UHF~algebras and therefore
$B$ embeds unitaly in $A$. Since $K_{0}(A)$ is unperforated,  
Theorem~\ref{thm_main} implies that $MU(A)$ exists, and so 
$B$ embeds unitaly into $MU(A)$. Consider 
the positive order unit preserving homomorphism 
$K_{0}(\iota) : K_{0}\left(MU(A)\right) \to K_{0}(A) \cong K_{0}(B)$. 
Then there is a unital $*$-homomorphism $MU(A) \to B$. Thus 
$MU(A) \cong B$ by Lemma~\ref{lem_uhf}\eqref{lem_uhf_it2}.
\end{proof}

Recall that a C*-algebra is called $K$-abelian if it is 
$KK$-equivalent to an abelian C*-algebra. The UCT class 
$\mathcal{N}$ is defined to be the family of 
all separable $K$-abelian C*-algebras 
\cite[Definition~2.4.5]{rs02}.

\begin{theorem}\label{thm_mu_kirchberg}
For any UHF~algebra $B$ there exists 
an uncountable family of pairwise non-isomorphic 
unital Kirchberg algebras in the UCT class $\mathcal{N}$ with 
a maximal UHF~subalgebra isomorphic to $B$.
\end{theorem}

\begin{proof}
Let $B$ be a UHF~algebra and consider 
the simple dimension group $K_{0}(B) + \alpha \mathbb{Z}$ 
where $\alpha$ is an arbitrary irrational number. According to 
\cite[Proposition~4.3.3(i)]{rs02}, there is 
a unital Kirchberg algebra $A(\alpha)$ in the UCT class 
$\mathcal{N}$ such that 
$K_{0}(A) \cong K_{0}(B) + \alpha \mathbb{Z}$ and 
$K_{1}(A) \cong \mathbb{Z}$. By 
Proposition~\ref{prop_list_K0lift}\eqref{prop_list_K0lift_it3}, 
$K_{0}(A(\alpha))$ has the $K_{0}$-lifting property for 
UHF~algebras, and hence by Theorem~\ref{thm_main}, 
$MU(A)$ exists. Similar to the first part of the proof of 
Theorem~\ref{thm_mu_af}, we get $MU(A(\alpha)) \cong B$.

By \cite[Theorem~8.4.1(iv)]{rs02} and similar to 
the second part of the proof of Theorem~\ref{thm_mu_af}, 
$\{ A(\alpha) : \alpha \in I \}$ is the desired family.
\end{proof}

\section{C*-algebraic Realization of the Rational Subgroup}\label{sec_rat_grp}

In this section we prove Theorems~\ref{thm_dim_k0_iso} and
\ref{thm_k0_dim_iso}, and Corollary~\ref{cor_mu_system}.
Recall that an ordered Abelian group $(G, G^{+})$ is said to be 
\emph{simple} if every nonzero $u$ in $G^{+}$ is 
an order unit (\cite[Definition~5.1.6]{rll00}). Recall that if 
$g \in G^{+}$ and $n\in\mathbb{N}$, then $n \vert g$ 
means that there is $x \in G^{+}$ such that $nx=g$.

\medskip

\begin{proposition}\label{prop_torfree_vert}
Let $G$ be a torsion-free Abelian group. Let $m, n \in \mathbb{N}$ 
be co-prime and $g \in G$.  If $mx = ng$ for some $x \in G$, 
then there exists $y \in G$ such that $my = g$. 
\end{proposition}

\begin{proof}
First we suppose that $G$ is countable. So let 
$G = \left\{ g_{1}, g_{2}, \ldots \right\} $. For any $j \geq 1$, 
set $G_{j} = \left\langle \{ g_{1}, \ldots, g_{j} \} \right\rangle$, 
the subgroup generated by $\{ g_{1}, \ldots, g_{n} \}$. Then 
$G_{j}$ is a finitely generated torsion-free subgroup of $G$ and 
using the fundamental theorem of finitely generated Abelian groups, 
it follows that $G_{j}$ is a free group. Thus there are some $n_{j}$ 
and a group isomorphism $\theta_{j} : G_{j} \to \mathbb{Z}^{n_{j}}$. 
Consider the direct limit 
\[ G_{1} 
\overset{j_{1}} \longrightarrow 
G_{2} 
\overset{j_{2}} \longrightarrow 
G_{3}
\overset{j_{3}} \longrightarrow 
\cdots \]
where $j_{n} : G_{j} \to G_{j+1}$ is the injection map. 
We have the following commutative diagram
\[
\xymatrix{
G_{1} \ar[r]^{j_{1}} \ar[d]^{\theta_{1}}_{\cong}
&G_{2} \ar[r]^{j_{2}} \ar[d]^{\theta_{2}}_{\cong}
&G_{3} \ar[r]^{j_{3}} \ar[d]^{\theta_{3}}_{\cong}
&\cdots \ar[r]
&G \ar[d]_{\cong}^{\theta}
\\
\mathbb{Z}^{n_{1}} \ar[r]^{\varphi_{1}} 
&\mathbb{Z}^{n_{2}} \ar[r]^{\varphi_{2}} 
&\mathbb{Z}^{n_{3}} \ar[r]^{\varphi_{3}} 
&\cdots \ar[r]
&H 
}
\]
where $\varphi_{n} = \theta_{n+1} \circ j_{n} \circ \theta_{n}^{-1}$. 
Each $\varphi_{n}$
is an injective group homomorphism and $H$ is the direct limit of 
the sequence $\left\{ \mathbb{Z}^{n_{j}}, \varphi_{j} \right\}_{j=1}^{\infty}$. 
By \cite[Propositions~6.2.5 and 6.2.6]{rll00}, $H = 
\bigcup_{j=1}^{\infty} \varphi^{j}(\mathbb{Z}^{n_{j}})$ where 
$\varphi^{j} : \mathbb{Z}^{n_{j}} \rightarrow H$ is 
the canonical injective group homomorphism. 

Now let $mx = ng$ and hence $m \theta(x) = n \theta(g)$.
Assume that $\theta(g) = \varphi^{k}(r)$ and 
$\theta(x) = \varphi^{l}(s)$ where $r \in \mathbb{Z}^{n_{k}}$, 
$s \in \mathbb{Z}^{n_{l}}$, and $k \leq l$. Since 
$\varphi^{l}\left(n \varphi_{l,k}(r)\right) = n \theta(g) = 
m \theta(x) = \varphi^{l}(ms)$ and $\varphi^{l}$ is injective, 
we have $n \varphi_{l,k}(r) = ms$ where 
$\varphi_{l,k} : \mathbb{Z}^{n_{k}} \to  \mathbb{Z}^{n_{l}}$ 
is the composition of $\varphi_{k}, \varphi_{k+1}, \ldots, \varphi_{l-1}$. 
Since \( \gcd(m, n) = 1 \), there exists $z \in \mathbb{Z}^{n_{l}}$ with 
$mz = \varphi_{l,k} (r)$ and hence $mw = \theta(g)$ where 
$w = \varphi^{l}(z)$. Letting $y = \theta^{-1}(w)$, 
we get $my = g$ as disared.

If $G$ is uncountable, a similar argument can be provided 
by taking $G_{F}$ the subgroup generated by $F$ where 
$F$ is a finite subset of $G$. Then $G$ is the direct limit of 
$G_{F}$'s and the rest of the argument works.
\end{proof}

\begin{corollary}\label{cor_dim_vert}
Let $(G, G^{+})$ be a dimension group. If  $g \in G^{+}$ 
and co-prime natural numbers $m$ and $n$ satisfy 
$m \vert n g$, then $m \vert g$. 
\end{corollary}

\begin{proof}
Let $mx = ng$ for some $x \in G^{+}$. By 
Proposition~\ref{prop_torfree_vert}, there is $y \in G$ 
such that $my = g$. Note that since $g \in G^{+}$, 
unperforation of $G$ implies that $y \in G^{+}$ and 
therefore $m \vert g$.  
\end{proof}

The notian of a ``rational subgroup’’ is defined for 
a simple dimension group in \cite{ghh17}. We define 
it for any ordered Abelian group.

\begin{definition}\label{def_rat_grp}
The \emph{rational subgroup} of an ordered Abelian group 
$(G, G^{+})$ with order unit $u$ is 
\[ \mathbb{Q}(G, u) := \left\{ g \in G : mg = qu \  \text{for some} \ 
m \in \mathbb{N} \ \text{and} \ q \in \mathbb{Z} \right\}. \]
\end{definition}

\begin{remark}
Let $(G, G^{+}, u)$ be an ordered Abelian group with 
order unit $u$. Then
$\mathbb{Q}(G, u) \subseteq G^{+} \cup -G^{+}$ and 
$\mathbb{Q}(G, u)$ is totally ordered group. In fact, 
if $g \in \mathbb{Q}(G, u)$, then there are 
$n \in \mathbb{N}$ and $p \in \mathbb{Z}$ 
such that $ng = pu$. If $p \geq 0$ then since 
$ng \geq 0$ and $G$ is unperforated, $g \geq 0$. 
If $p < 0$ then $n (-g) \geq 0$ and hence $g < 0$.
To see that $\mathbb{Q}(G, u)$ is totally ordered group, 
let $g$ and $h$ be in $\mathbb{Q}(G, u)$. 
Consider integers $n, m \in \mathbb{N}$ 
and $p, q$ in $\mathbb{Z}$ such that 
$ng = pu$ and $mh = qu$. We may assume that $qn \geq pm$.  
If $g, h \in G^{+}$ then $pm \geq 0$ implies that
$pmg \leq qng = pqu = pmh$ and therefore $g \leq h$ (since $G$ 
is unperforated). Also if $g, h \in -G^{+}$  
then $-pmg \geq -pmh$ and therefore $g \geq h$ 
(since $-pm \geq 0$). The cases $g \geq 0 \geq h$ and 
$h \geq 0 \geq g$ imply $g \geq h$ and $h \geq g$, respectively.
\end{remark}

\begin{lemma}\label{lem_q(q)}
Let $\alpha$ be an irrational number and $H$ be a subgroup of 
$\mathbb{Q}$ such that $1 \in H$. Let $v$ is a positive element 
in the additive subgroup $H + \alpha \mathbb{Z}$ of $\mathbb{R}$ 
and $v = k + \alpha z$ for some (necessarily unique) 
$k \in H \backslash \{ 0 \}$ and $z \in \mathbb{Z}$.  
Then 
\[ \mathbb{Q}(H + \alpha \mathbb{Z}, v) = \left\{ h + \alpha \frac{hz}{k} : 
h \in H \ \text{and} \ \frac{hz}{k} \in \mathbb{Z} \right\}. \]
In particularly, $\mathbb{Q}(H + \alpha \mathbb{Z}, 1) = H$.
\end{lemma}

\begin{proof}
First note that $H + \alpha \mathbb{Z}$ is a simple dimension group. 
Let $g = h + \alpha w$ be in 
$\mathbb{Q}(H + \alpha \mathbb{Z}, v)$ where $h \in H$ and 
$w \in \mathbb{Z}$. There are $m \in \mathbb{N}$ and 
$q \in \mathbb{Z}$ such that $mg = qv$. 
So $mh + \alpha mw = qk + \alpha qz$ and hence
$mh = qk$ and $mw = qz$. Hence $w = hz/k$. 

Conversely, let $h + (\alpha hz/k)$ be in the right hand set in the statement. 
Take $m \in \mathbb{N}$ and $q \in \mathbb{Z}$ with $h/k = q/m$. Then 
$m(h + (\alpha hz/k)) = qv$ and therefore $h + (\alpha hz/k)$ belongs to 
$\mathbb{Q}(H + \alpha \mathbb{Z}, v)$. 

For the last part of the statement, let $v = 1$. Then $k = 1$ and 
$z = 0$. Hence $\mathbb{Q}(H + \alpha \mathbb{Z}, 1) = H$.
\end{proof}

Recall from Subsection~\ref{subsec_propd} the definition of 
the supernatural number $N(G, u)$ associated to 
an ordered Abelian group with order unit $(G, G^{+}, u)$, 
and recall from Subsection~\ref{subsec_ord_grp} 
the subgroup $Q(N)$ of $\mathbb{Q}$ associated to 
a supernatural number $N$. Note that $1 \in Q(N)$. 
Then the following result is immediate from the preceding lemma.

\begin{corollary}
Let $(G, u)$ be a dimension group with order unit $u$ and 
$\alpha$ be an irrational number. Then 
$\mathbb{Q}\left(\mathit{Q}(N(G, u)) + \alpha \mathbb{Z}, 1\right) 
= \mathit{Q}(N(G, u))$.
\end{corollary}

Let $(G, G^{+}, u)$ be 
an ordered Abelian group with a distinguished order unit and 
consider the supernatural number 
$N(G, u) = \{ n_{j} \}_{j = 1}^{\infty}$. Note that 
an arbitrary element of $\mathit{Q}\left(N(G, u)\right)$ is written as 
$\sum_{j=1}^{k} \alpha_{j}/p_{j}^{m_{j}}$ where 
$m_{j} = n_{j}$ if $n_{j} \neq \infty$ and $m_{j}$ is an
arbitrary nonnegative integer if $n_{j} = \infty$. We define 
a homomorphism  
\begin{align*}
\theta 
&: \mathit{Q}(N(G, u)) \rightarrow \mathbb{Q}(G, u) 
\\
&\hspace{0.93 cm} \sum\limits_{j=1}^{k} \frac{\alpha_{j}}{p_{j}^{m_{j}}} 
\mapsto \sum\limits_{j=1}^{k} \alpha_{j} x(j, m_{j})
\end{align*}
where $x(j, m_{j})$ is the unique element of 
$G^{+}$ with $p_{j}^{m_{j}} x(j, m_{j}) = u$ for all 
$1 \leq j \leq k$. 

\begin{theorem}\label{thm_grp_iso_2}
Let $(G, u)$ be a dimension group with order unit $u$. Then 
the map 
$\theta : (\mathit{Q}(N(G, u)), 1) \rightarrow \mathbb{Q}(G, u)$ 
defined above, is an isomorphism of dimension groups with order unit. 
\end{theorem}

\begin{proof}
First we show that $\theta$ is well defined. For integers 
$\alpha_{1}, \ldots, \alpha_{k}$ and nonnegative integers 
$m_{1}, \ldots, m_{k}$ with 
$p = p_{1}^{m_{1}} \cdots p_{k}^{m_{k}}$, since
\begin{equation}\label{equ_main}
p \sum\limits_{j=1}^{k} \alpha_{j} x(j, m_{j}) = 
p \Big{(} \sum\limits_{j=1}^{k} \frac{\alpha_{j}}{p_{j}^{m_{j}}} \Big{)} u,
\end{equation}
we see that $\sum_{j=1}^{k} \alpha_{j} x(j, m_{j})$ belongs to 
$\mathbb{Q}(G, u)$. Also for integers 
$\beta_{1}, \ldots, \beta_{k}$ and nonnegative integers 
$n_{1}, \ldots, n_{k}$, if
$\sum_{j=1}^{k} \alpha_{j}/p_{j}^{m_{j}} 
= \sum_{j=1}^{k} \beta_{j}/p_{j}^{n_{j}}$
then 
$\sum_{j=1}^{k} \big{(} (\alpha_{j}/p_{j}^{m_{j}}) - 
(\beta_{j}/p_{j}^{n_{j}}) \big{)} = 0$. 
Hence  
$q \Big{(} \sum_{j=1}^{k} \big{(} (\alpha_{j}/p_{j}^{m_{j}}) - 
(\beta_{j}/p_{j}^{n_{j}}) \big{)} \Big{)} u 
= 0$ 
where 
$q = p_{1}^{m_{1} + n_{1}} \cdots p_{k}^{m_{k} + n_{k}}$.
Using  $p_{j}^{m_{j}} x(j, m_{j}) = u$ for all $1 \leq j \leq k$,
it follows that
$\sum_{j=1}^{k} \alpha_{j} x(j, m_{j}) = 
\sum_{j=1}^{k} \beta_{j} x(j, n_{j})$. Thus $\theta$ is well defined.

It follows from Equation~\eqref{equ_main} and that 
$\mathbb{Q}(G, u)$ is unperforated, $\theta$ and $\theta^{-1}$ 
are positive. Also $\theta$ is order unit preserving, because 
\[ \theta(1) = \theta\Bigg{(}\frac{p_{j}^{m_{j}}}{p_{j}^{m_{j}}}\Bigg{)} 
= p_{j}^{m_{j}} x(j, m_{j}) = u \] 
where $1 \leq j \leq k$. Equation~\eqref{equ_main} and unperforation 
imply that $\theta$ is injective. For surjectivity, 
first we prove the following claim:

\medskip

\textbf{Claim}. Let $g$ be in $\mathbb{Q}(G, u)$ with 
$mg = qu$ for some $m \in \mathbb{N}$ and 
$q \in \mathbb{Z}$.  Let $p = p_{1}^{m_{1}} \cdots p_{k}^{m_{k}}$ 
be an integer where $p_{j}$'s are distinct prime numbers and $m_{j}$'s 
are natural numbers. Then $g = \theta(x)$ for some 
$x = \sum_{j=1}^{k} \big{(}\alpha_{j}/p_{j}^{m_{j}}\big{)}$
in $\mathit{Q}(N(G, u))$ if and only if $pq/m$ is an integer.

\medskip

To prove this claim, first assume that $g = \theta(x)$ for some 
$x = \sum_{j=1}^{k} \big{(}\alpha_{j}/p_{j}^{m_{j}}\big{)}$
in $\mathit{Q}(N(G, u))$. Since
\begin{align*}
pqu 
&= pm \theta(x) = pm \sum\limits_{j=1}^{k} \alpha_{j} x(j, m_{j}) \\
&= m \sum\limits_{j=1}^{k} \alpha_{j} p x(j, m_{j}) \\
&= m \Bigg{(} \sum\limits_{j=1}^{k} \Big{(} \alpha_{j} \prod_{\substack{1 \leq i \leq k \\ i \neq j}} 
p_{i}^{m_{i}} \Big{)} \Bigg{)} u,
\end{align*} 
we have $pq/m \in \mathbb{Z}$. Conversely, let $pq/m \in \mathbb{Z}$. 
Choose integers $\beta_{1}, \ldots, \beta_{k}$ with
$\sum_{j=1}^{k} \big{(} \beta_{j} \prod_{\substack{1 \leq i \leq k \\ i \neq j}} p_{i}^{m_{i}} \big{)} = 1$. 
Take $x = \sum_{j=1}^{k} \big{(}\alpha_{j}/p_{j}^{m_{j}}\big{)}$ 
where $\alpha_{j} = \beta_{j} pq /m$ for $1 \leq j \leq k$. Since 
$m \sum_{j=1}^{k} \big{(} \alpha_{j} 
\prod_{\substack{1 \leq i \leq k \\ i \neq j}} p_{i}^{m_{i}} \big{)} = pq$ 
and
\begin{align*}
mpg 
&= pqu = m \Bigg{(} \sum\limits_{j=1}^{k} \Big{(} \alpha_{j} 
\prod\limits_{\substack{1 \leq i \leq k \\ i \neq j}} p_{i}^{m_{i}} \Big{)} \Bigg{)} u \\ 
&= m \Big{(} \sum\limits_{j=1}^{k} \alpha_{j} p x(j, m_{j}) \Big{)} = mp \theta(x),
\end{align*}
we get that $g = \theta(x)$, as $G$ is unperforated. This finishes the proof of the claim.

\medskip

Now by this claim, we show that the map $\theta$ is surjective. 
Let for some $g \in \mathbb{Q}(G, u)^{+}$ we have 
$mg = qu$ where 
$m = p_{1}^{r_{1}} \cdots p_{k}^{r_{k}}$ and 
$q = p_{1}^{s_{1}} \cdots p_{k}^{s_{k}}$ are 
prime factorizations. Consider 
the supernatural number 
$N(G, u) = \{ n_{j} \}_{j \in \mathbb{N}}$ and 
$p = p_{1}^{m_{1}} \cdots p_{k}^{m_{k}}$ where 
$m_{j} = n_{j}$ if $n_{j} \neq \infty$ and $m_{j} = r_{j} - s_{j}$ 
if $n_{j} = \infty$. Now if $n_{j} = \infty$ then
$p_{j}^{s_{j} + m_{j} - r_{j}} \in \mathbb{Z}$ for all $1 \leq j \leq k$ 
with $m_{j} = r_{j} - s_{j}$. Now let $n_{j} \neq \infty$ for some 
$1 \leq j \leq k$, we show that $r_{j} \leq s_{j} + n_{j}$. Contrary 
suppose that $r_{j} > s_{j} + n_{j}$. Since 
$u = p_{j}^{n_{j}} x(j, n_{j})$,
we have
\begin{equation*}
p_{j}^{r_{j}} \prod_{\substack{1 \leq i \leq k \\ i \neq j}} p_{i}^{r_{i}} g 
= mg = qu = p_{j}^{s_{j} + n_{j}} 
\prod_{\substack{1 \leq i \leq k \\ i \neq j}} p_{i}^{s_{i}} x(j, n_{j})
\end{equation*}
and hence by unperforation
\begin{equation*}
p_{j}^{r_{j} - s_{j} - n_{j}} \prod_{\substack{1 \leq i \leq k \\ i \neq j}} 
p_{i}^{r_{i}} g = \prod_{\substack{1 \leq i \leq k \\ i \neq j}} 
p_{i}^{s_{i}} x(j, n_{j}).
\end{equation*}
Since $p_{j}^{r_{j} - s_{j} - n_{j}}$ and 
$\prod_{\substack{1 \leq i \leq k \\ i \neq j}}  p_{i}^{s_{i}}$ are 
relatively prime, Corollary~\ref{cor_dim_vert} implies that 
$p_{j}^{r_{j} - s_{j} - n_{j}}$ devides $x(j, n_{j})$. 
Hence there is $h \in \mathbb{Q}(G, u)^{+}$ such that
$p_{j}^{r_{j} - s_{j} - n_{j}} h = x(j, n_{j})$ 
and hence 
$p_{j}^{r_{j} - s_{j}} h = p_{j}^{n_{j}} x(j, n_{j}) = u$, 
by the definition of $n_{j}$. Thus we get 
$r_{j} - s_{j} \leq n_{j}$ but it is a contradiction. Then 
$r_{j} \leq s_{j} + n_{j}$, and hence 
$p_{j}^{s_{j} + m_{j} - r_{j}} \in \mathbb{Z}$ for all 
$1 \leq j \leq m$.

\medskip

Finally we see that 
$pq/m = \prod_{1 \leq i \leq k} p_{j}^{s_{j} + m_{j} - r_{j}}$ 
belongs to $\mathbb{Z}$. By the claim, $g$ belongs to 
the range of the map $\theta$ and thus this map is surjective.  
\end{proof}

We are ready to give the following proofs.

\begin{proof_thmxx}
Let $N = N(G, u)$ be the supernatural number of $(G, u)$ 
as in Subsection~\ref{subsec_propd}. Consider 
the dimension group $Q(N)$. By \cite[Proposition~7.2.8]{rll00},
there is an AF~algebra $A$ such that 
$(K_{0}(A), [1_{A}]_{0}) \cong (G, u)$. 
By \cite[Proposition~4.3.4]{rs02}, 
there is a unital Kirchberg algebra $B$ in the UCT calss 
$\mathcal{N}$ such that $(K_{0}(B), [1_{B}]_{0}) \cong (G, u)$ 
and $K_{1}(B) = \mathbb{Z}$. By Theorem~\ref{thm_main}, 
C*-algebras $A$ and $B$ have maximal UHF~subalgebras $MU(A)$ 
and \( MU(B) \), respectively. By Lemma~\ref{lem_weakunp_propd} 
and Theorem~\ref{thm_max_propd}, the set $\Sigma_{A}$ of 
supernatural numbers $M$ with $M \vert [1_{A}]_{0}$, has 
the maximum element. Since $(K_{0}(A), [1_{A}]_{0}) \cong (G, u)$, 
\( N \) is the maximum element of $\Sigma_{A}$ and by the proof of 
Theorem~\ref{thm_main}, \( K_{0}(MU(A)) \cong Q(N) \). Therefore, 
\[ K_{0}(MU(A)) \cong K_{0}(MU(B)) \cong Q(N) \cong \mathbb{Q}(G, u), \]
and this finishes the proof. \hfill
\end{proof_thmxx}

Combining Theorem~\ref{thm_main} and Theorem~\ref{thm_grp_iso_2}, 
we are able to give the proof of 
Theorem~\ref{thm_k0_dim_iso}: 

\begin{proof_thmxxx}
Since $K_{0}(A)$ is a dimension group, by Theorem~\ref{thm_grp_iso_2}
$\mathbb{Q}(K_{0}(A), [1]_{0})$ is isomorphic to 
$\mathit{Q}(N(K_{0}(A), [1]_{0}))$. Also according to 
the proof of Theorem~\ref{thm_main}, 
$\mathit{Q}(N(K_{0}(A), [1]_{0})$ is isomorphic to 
$K_{0}(MU(A))$ and therefore, 
\[ (K_{0}(MU(A)), [1]_{0}) \cong \mathit{Q}(N(K_{0}(A), [1]_{0})) \cong 
\mathbb{Q}(K_{0}(A), [1]_{0}), \] 
as desired. \hfill
\end{proof_thmxxx}

\medskip

\begin{example}\label{exa_mu_part2}
$MU(\mathcal{O}_{\infty}) \cong \mathbb{C}$. For this, 
since $K_{0}(\mathcal{O}_{\infty}) = \mathbb{Z}$
(\cite{bl86}), by 
Proposition~\ref{prop_list_K0lift}\eqref{prop_list_K0lift_it3} 
and
Theorem~\ref{thm_main}, $\mathcal{O}_{\infty}$ has 
a maximal UHF~subalgebra. Let 
$MU(\mathcal{O}_{\infty}) \cong \mathrm{M}_{N}$. By 
Theorem~\ref{thm_k0_dim_iso}, 
\begin{align*}
(Q(N), 1) 
&\cong (K_{0}(MU(\mathcal{O}_{\infty})), [1_{\mathcal{O}_{\infty}}]_{0}) \\
&\cong \mathbb{Q}(K_{0}(\mathcal{O}_{\infty}), [1_{\mathcal{O}_{\infty}}]_{0}) \\
&\cong \mathbb{Q}(\mathbb{Z}, 1) \\
&\cong Q(N(\mathbb{Z}, 1)),
\end{align*}
and by \cite[Proposition~7.4.3(ii)]{rll00}, 
$N = \{ 0, 0, \cdots \}$. Therefore 
$MU(\mathcal{O}_{\infty}) \cong \mathbb{C}$. 
\end{example}

\begin{proof_corx}
Let $(X, T)$ be a Cantor minimal system and consider 
the C*-algebra crossed product  
$A = C(X) \rtimes_{T} \mathbb{Z}$. In \cite{pu89} it is 
shown that the group $K^{0}(X, T)$ is order isomorphic to 
the group $K_{0}(A)$. Therefore, 
\[ \mathbb{Q}(K_{0}(A), [1]_{0}) \cong 
\mathbb{Q}(K^{0}(X, T), [1_{X}]). \] 
The C*-algebra $A$ is unital separable simple tracial rank zero. 
Also according to \cite[Theorem~4.1]{pu89}, $K^{0}(X, T)$ is 
a simple, acyclic (i.e, $G \ncong \mathbb{Z}$) 
dimension group with (canonical) distinguished order unit $1$. 
Now Theorem~\ref{thm_main} imlies that $A$ has  
a maximal UHF~subalgebra $MU(A)$. By 
Theorem~\ref{thm_k0_dim_iso}, 
$\mathbb{Q}(K_{0}(A), [1]_{0})) \cong (K_{0}(MU(A)), [1]_{0})$. 
Also \cite[Proposition~3.31]{ho16} implies that 
$\mathbb{Q}(K^{0}(X, T), 1_{X}) \cong (K^{0}(Y, S), 1)$. 
Thus we have
\[ (K_{0}(MU(A)), 1) \cong \mathbb{Q}(K_{0}(A), [1]_{0})) \cong 
\mathbb{Q}(K^{0}(X, T), 1) \cong (K^{0}(Y, S), 1), \]
as desired. \hfill
\end{proof_corx}

\section{Maximal UHF~subalgebras of AF~algebras}
\label{sec_mu_bra}

In this section, we give another method to prove 
Theorem~\ref{thm_main} for unital AF~algebras in which 
maximal UHF~subalgebras are obtained by a combinatorial method 
using Bratteli diagrams. In practice, given an AF~algebra \( A \), 
first we draw its Bratteli diagram \( \mathcal{B}(A) \). Second,
 we draw a Bratteli diagram \( \mathcal{O}(\mathcal{B}(A)) \) 
 associated to \( \mathcal{B}(A) \) as described in 
 \cite[Definition~4.11]{aeg21} which has only one vertex at each level. 
 Finally, the UHF~algebra whose Bratteli diagram is 
 \( \mathcal{O}(\mathcal{B}(A)) \) is (up to isomorphism) 
 the desired maximal UHF~subalgebra of \( A \).

\medskip

A Bratteli diagram can be defined in two equivalent ways: 
using directed graphs \cite{br72, hps92} and using 
the matrix language \cite{aeg15}. We follow the first one here.
Let us recall the definition of a Bratteli diagram and 
a premorphism.

\begin{definition}[\cite{aeg21}, Definition~2.1]\label{def_brat_graph}
A \emph{Bratteli diagram} consists of a vertex set $V$ and 
an edge set $E$ satisfying the following conditions. We have 
a decomposition of $V$ as a disjoint union 
$V_0 \cup V_1 \cup \cdots$, where each $V_n$ is finite and 
nonempty and $V_0$ has exactly one element, $v_0$. Similarly, 
$E$ decomposes as a disjoint union 
$E_{1} \cup E_{2} \cup \cdots$,
where each $E_n$ is finite and nonempty. Moreover, we have 
maps $r, s : E \to V$ such that $r(E_n) \subseteq V_n$ and 
$s(E_n) \subseteq V_{n-1}, n=1,2,3, \ldots$ ($r =$ range, 
$s =$ source). We also assume that $s^{-1}(v)$ is nonempty 
for all $v$ in $V$ and $r^{-1}(v)$ is nonempty for all $v$ in 
$V \backslash V_0$. We denote the matrix associated with 
each edge set $E_{n}$ by $\mathrm{M}(E_{n})$ and call 
the \emph{multiplicity matrix} of $E_{n}$.

Note that each $\mathrm{M}(E_{n})$ is an embedding matrix 
in the sense that for each $j$ there is an $i$ such that 
the $ij-$th entry of $\mathrm{M}(E_{n})$ is nonzero.
\end{definition}

\begin{definition}[\cite{aeg21}, Definition~2.5]\label{def_pmor_graph}
Let $B = (V, E)$ and $C = (W, S)$ be Bratteli diagrams. By 
a \emph{premorphism} $f : B \rightarrow C$, we mean a pair 
$(F, (f_n)_{n = 0}^{\infty})$ where $(f_n)_{n = 0}^{\infty}$ 
is a cofinal (i.e., unbounded) sequence of positive integers with 
$f_{0} = 0 \leq f_{1} \leq f_2 \leq \cdots$, $F$ consists of 
a disjoint union $F_0 \cup F_1 \cup F_2 \cup \cdots$, together 
with a pair of range and source maps 
$r :F \rightarrow W, s : F \rightarrow V$ such that 
the following hold:
\begin{enumerate}
\item
each $F_n$ is a nonempty finite set, 
$s(F_n) \subseteq V_n, r(F_n) \subseteq W_{f_n} , F_0$ 
is a singleton, $s^{-1} \{ v \}$, is nonempty for all $v$ in $V$, 
and $r^{-1} \{w \}$ is nonempty for all $w$ in $W$;

\item
the diagram of $f:B \rightarrow C$,
\[ \xymatrix{
\ar @{} [dr] |{}
V_{0} \ar[d]_{F_{0}} \ar^{E_{0}}[r] 
& 
V_{1} \ar[d]_{F_{1}} \ar^{E_{1}}[r]  
& 
V_{2} \ar[d]_{F_{2}} \ar^{E_{2}}[r]  
& 
\cdots 
\\
W_{F_{0}} \ar_{S_{f_{0}, f_{1}}}[r] 
& 
W_{F_{1}} \ar_{S_{f_{1}, f_{2}}}[r]  
& 
W_{F_{2}} \ar_{S_{f_{2}, f_{3}}}[r]  
& 
\cdots 
} \]
commutes. The commutativity of the diagram of $f$ means that 
$E_{n+1} \circ F_{n+1} \cong F_n \circ S_{f_n, f_{n+1}}$ 
for all $n \geq 0$, i.e., there is a bijective map from 
$E_{n+1} \circ F_{n+1}$ to $F_n \circ S_{f_n, f_{n+1}}$ 
preserving the respective source and range maps.
\end{enumerate}
\end{definition}

We recall the Bratteli diagram \( \mathcal{B}(A) \) of 
a unital AF~algebra \( A \) \cite{br72}. Let \( A \) be 
the inductive limit of a sequence 
\( \{ (A_{n}, \varphi_{n}) \}_{n=0}^{\infty} \) 
where \( A_{0} \cong \mathbb{C} \), each \( A_{n} \) is 
a finite dimensional C*-algebra, and each 
\( \varphi_{n} : A_{n} \to A_{n+1} \)
is a unital \( * \)-homomorphism. The Bratteli diagram 
\( \mathcal{B}(A) = (V, E) \) of \( A \) (depending on 
\( A_{n} \)'s and \( \varphi_{n} \)'s)
has the vertex set \( V = \bigcup_{n=0}^{\infty} V_{n} \) 
where \( V_{0} \) is a singleton and \( \# V_{n} \) 
equals the number of full matrix algebra summands whose 
direct sum is isomorphic to \( A_{n} \). Each edge set \( E_{n} \) 
is obtained from the multiplicity matrix of \( \varphi_{n} \)
according to \cite[Theorem~2.1]{aeg15}. Note that, 
though the Bratteli diagram of \( A \) is not unique 
(as it depend on the inductive system), any two Bratteli diagrams of 
\( A \) are equivalent \cite{br72, aeg15}.  

The equivalence of first two parts following result is 
the special case of some results of \cite{aeg15} for 
unital AF~algebras (see Section~3 and the proof of 
Theorem~4.1 in \cite{aeg15}). 

\begin{proposition}\label{prop_af_lift}
Let $A$ and $B$ be unital AF~algebras and $\mathcal{B}(A)$ 
and $\mathcal{B}(B)$ be Bratteli diagrams for $A$ and $B$, 
respectively. Then the following statements are equivalent:
\begin{enumerate}
\item\label{prop_af_lift_it1}
there is a premorphism 
$f: \mathcal{B}(A) \rightarrow \mathcal{B}(B)$,

\item\label{prop_af_lift_it2}
there exists a unital $*$-homomorphism $\varphi : A \to B$,

\item\label{prop_af_lift_it3}
there is a positive group homomorphism 
$\alpha : K_{0}(A) \to K_{0}(B)$ such that 
\( \alpha([1_{A}]_{0}) = [1_{B}]_{0} \) and
$K_{0}(\varphi) = \alpha$.
\end{enumerate}
\end{proposition}

\begin{proof}
The equivalence of \eqref{prop_af_lift_it1} and 
\eqref{prop_af_lift_it2} is given in \cite{aeg15}.  
The equivalence of 
\eqref{prop_af_lift_it2} and \eqref{prop_af_lift_it3} follows from 
Paragraph~3.2.2 and Exercise~7.7 of \cite{rll00}.
\end{proof}

\begin{theorem}\label{thm_main_af}
Let \( A \) be a unital AF~algebra. Then a maximal UHF~subalgebra 
$MU(A)$ of $A$ exists. Moreover, for any 
UHF unital C*-subalgebra $D$ of $A$, there exists 
a unital embedding $\phi: D \rightarrow MU(A)$ with 
$\iota_{MU(A)} \circ \phi \approx_{a. u.} \iota_{D}$ where
\( \iota_{D} \) denotes the injection map from \( D \) to \( A \).
\end{theorem}

\begin{proof}
\emph{Existence}: 
There is an inductive limit
\[ \mathbb{C}1_{A} 
\overset{\iota_{\mathbb{C}1_{A}}} \longrightarrow 
A_{1} 
\overset{\iota_{A_{1}}} \longrightarrow 
A_{2} 
\overset{\iota_{A_{2}}} \longrightarrow 
\cdots 
\longrightarrow 
A, \]
where $A_{n}$ is a finite dimensional C*-subalgebra of $A$ 
and $\iota_{A_n}$ is the inclusion for all $n \geq 1$. 
We consider the Bratteli diagram $\mathcal{B}(A)$ of $A$
as described before Proposition~\ref{prop_af_lift}. 
Consider the odometer $\mathcal{O}(\mathcal{B}(A)) = (W, R)$ 
of type $(r_{n})^{\infty}_{n=1}$, and the premorphism 
$f_{\mathcal{B}(A)} : \mathcal{O}(\mathcal{B}(A)) \to \mathcal{B}(A)$
as in \cite[Definition~4.11]{aeg21}. 
To recall, Let  $\mathrm{M}(E_{n})$ denote the multiplicity matrix
of $E_{n}$. Then $E_{0,n}$ defined by 
$E_{1}\circ E_{2}\circ\cdots\circ E_{n}$ (the edge set from $V_{0}$ to 
$V_{n}$) is the set of towers at  level $n$, and the column matrix
\[
\mathrm{M}(E_{0,n})= \mathrm{M}(E_{n}) \cdots \mathrm{M}(E_{n-1})\mathrm{M}(E_{1})
=\left(
\begin{smallmatrix}
h_{n,1}\\
h_{n,2}\\
\vdots\\
h_{n,k_{n}}
\end{smallmatrix}
\right),
\]
where the $h_{n,i}$ are non-zero positive integers and 
$k_{n}= \# V_{n}$, consists of the heights of these towers.
We set $h_{n}=\mathrm{gcd}(h_{n,1},h_{n,2}, \ldots, h_{n,k_{n}})$.
Note that $1=h_{0}\mid h_{1}\mid h_{2}\cdots$
and so the definition of $r_{n}={h_{n}}/{h_{n-1}}$ makes sense.

Let $B$ be the UHF~algebra
whose Bratteli diagram is \( \mathcal{O}(\mathcal{B}(A)) \), 
more precisely, \( B \) is the inductive limit of 
the following inductive sequence  
\[ \mathbb{C} 
\overset{\psi_{0}} \longrightarrow 
B_{1} 
\overset{\psi_{1}} \longrightarrow 
B_{2} 
\overset{\psi_{2}} \longrightarrow 
\cdots, \]
where $B_n = M_{r_1 r_2 \ldots r_n}$ for $n \ge 1$ and
$\psi_n : B_n \to B_{n+1}$ is the $*$-homomorphism 
defined by 
$\psi_{n}(a) = \text{diag}(a, \ldots, a)$ (with $r_{n+1}$ 
copies of $a$) for \( n \geq 0 \) where \( B_{0} = \mathbb{C} \). 
Since $\mathcal{B}(B) = \mathcal{O}(\mathcal{B}(A))$, by
Proposition~\ref{prop_af_lift} there exists 
a unital $*$-homomorphism $\varphi : B \to A$. Define 
$MU(A) = \varphi(B)$.

\medskip

\emph{Maximality}: Let $D$ be a UHF~subalgebra of $A$ 
with \( 1_{A} \in D \) and consider the Bratteli diagram $\mathcal{B}(D)$ 
and the premorphism $g : \mathcal{B}(D) \to \mathcal{B}(A)$ 
associated to the unital \( * \)-homomorphism 
\( \iota_{D} : D \to A \) as in \cite[Definition~3.3]{aeg15}. By 
the proof of \cite[Theorem~4.12]{aeg21}, there exists 
a premorphism 
$h : \mathcal{B}(D) \to \mathcal{B}(MU(A))$ such that 
$f_{\mathcal{B}(A)} \circ h = g$. By 
Proposition~\ref{prop_af_lift}, there is 
a unital $*$-homomorphism $\phi : D \to MU(A)$ and by 
\cite[Lemma~5.4]{aeg15},
$\iota_{MU(A)} \circ \phi \approx_{a.u.} \iota_{D}$.
\end{proof}

\begin{example}\label{exa_mu_part3}
Let $A$ be the inductive limit of the following sequence
\[ \mathbb{C}1 
\overset{\psi_{0}} \longrightarrow 
\mathbb{C}1 \oplus \mathbb{C}1 
\overset{\psi_{1}} \longrightarrow 
\mathrm{M}_{3^{2}} \oplus \mathrm{M}_{3^{2}} 
\overset{\psi_{2}} \longrightarrow 
\cdots, \]
where 
$\psi_{n}(x \oplus y) = 
\mathrm{diag}(x, x, y) \oplus \mathrm{diag}(x, y, y)$ 
for all $x, y \in \mathrm{M}_{3^{n}}$. Then by the proof of
Theorem~\ref{thm_main_af}, 
$MU(A) \cong \mathrm{M}_{3^{\infty}}$. In fact, the Bratteli diagram
of \( A \) is on the right in Fiqure~\ref{fig}. The diagram 
\( \mathcal{O}(\mathcal{B}(A)) \) and the premorphism 
\( f_{\mathcal{B}(A)} : \mathcal{O}(\mathcal{B}(A)) \to 
\mathcal{B}(A) \) described in the proof of 
Theorem~\ref{thm_main_af}, are depicted in Figure~\ref{fig}. 
Note that \( r_{n} = 3 \) since \( h_{n} = 3^{n-1} \) and therefore 
\( r_{n} = h_{n} / h_{n-1} = 3^{n-1} / 3^{n-2} = 3 \) 
for all \( n \geq 2 \) and \( r_{1} = h_{0} = 1 \). 
Note that \( A \) is not a UHF~algebra since otherwise 
it implies that \( MU(A) \cong A \) and hence 
there is a premorphism from \( \mathcal{B}(A) \) to 
\( \mathcal{O}(\mathcal{B}(A)) \), by 
Proposition~\ref{prop_af_lift}. However, looking at Figure~\ref{fig}, 
by inspection there is no premorphism from the right Bratteli diagram to
the left. 

\begin{figure}[h]
\begin{center}
\begin{tikzpicture}[scale=1.1]

\filldraw (3,10) circle [radius=0.1];
\filldraw (3,7) circle [radius=0.1];
\filldraw (3,4) circle [radius=0.1];
\filldraw (3,1) circle [radius=0.1];

\draw (3,9.87)--(3,7.14);

\draw (2.9,6.87)--(2.9,4.14);
\draw (3,6.87)--(3,4.14);
\draw (3.1,6.87)--(3.1,4.14);

\draw (2.9,3.87)--(2.9,1.14);
\draw (3,3.87)--(3,1.14);
\draw (3.1,3.87)--(3.1,1.14);

\filldraw (7,10) circle [radius=0.1];
 \filldraw (6,7) circle [radius=0.1];
\filldraw (8,7) circle [radius=0.1];
 \filldraw (6,4) circle [radius=0.1];
\filldraw (8,4) circle [radius=0.1];
 \filldraw (6,1) circle [radius=0.1];
\filldraw (8,1) circle [radius=0.1];


\draw (6.94,9.9)--(6.01,7.13);
\draw (7.085,9.915)--(7.97,7.12);

\draw (5.97,6.87)--(5.97,4.14);
\draw (5.97,3.87)--(5.97,1.14);
\draw (6.06,6.87)--(6.06,4.14);
\draw (6.06,3.87)--(6.06,1.14);

\draw (7.98,6.87)--(7.98,4.14);
\draw (7.98,3.87)--(7.98,1.14);
\draw (8.06,6.87)--(8.06,4.14);
\draw (8.06,3.87)--(8.06,1.14);

\draw (7.9,6.9)--(6.08,4.1);
\draw (7.9,3.9)--(6.08,1.1);
\draw (6.1,6.9)--(7.95,4.15);
\draw (6.1,3.9)--(7.95,1.15);

\draw[->, thick] (3.1,10.1) [out=10,in=170] to (6.9,10.1);

\draw[->, thick] (3.1,7.1) [out=10,in=170] to (5.9,7.1);
\draw[->, thick] (3.1,7.15) [out=30,in=150] to (7.9,7.1);

\draw[->, thick] (3.15,4.1) [out=10,in=170] to (5.9,4.1);
\draw[->, thick] (3.15,4.15) [out=30,in=150] to (7.9,4.1);

\draw[->, thick] (3.15,1.1) [out=10,in=170] to (5.9,1.1);
\draw[->, thick] (3.15,1.15) [out=30,in=150] to (7.9,1.1);
\node[very thick] at (3,0.2) {\vdots};
\node[very thick] at (7,0.2) {\vdots};
\node at (3,10.7) {$\mathcal{O}(\mathcal{B}(A))$};
\node at (7,10.7) {$\mathcal{B}(A)$};
\node at (5,10.9) {$f_{\mathcal{B}(A)}$};
\draw[->] (4,10.7) to (6,10.7);
\end{tikzpicture}
\end{center}
\caption{}\label{fig}
\end{figure}
\end{example}

In the following remark we compare the notion of a maximal UHF~subalgebra 
in the sense of Definition~\ref{def_mu} and the same notion with respect to inclusion.

\begin{remark}\label{rem_main}
\begin{enumerate}
\item\label{rem_main_it1}
Let $A$ be a separable unital C*-algebra and
\[ \mathcal{U} = \left\{ D : D \ \text{is a UHF~C*-subalgebra of \( A \) and} 
\ 1_{D} = 1_{A} \right\}. \] Then by the Zorn's lemma and the fact that 
every separable unital C*-algebra which is locally UHF~algebra is indeed 
a UHF~algebra, the set \( \mathcal{U} \) has at least one maximal element with 
respect to inclusion. If a maximal UHF~subalgebra $MU(A)$ of $A$ as in  
Definition~\ref{def_mu} exists then $MU(A)$ is isomorphic to
a maximal element of $\mathcal{U}$. Indeed, 
the subset \( \mathcal{U}^{\prime} \) of \( \mathcal{U} \) consisting of 
elements \( D \in \mathcal{U} \) with \( MU(A) \subseteq D \), has 
a maximal element, say \( B \). Since $MU(A) \subseteq B$ and
$B$ is embedded in $MU(A)$, by Lemma~\ref{lem_uhf}, 
$MU(A) \cong B$. Also $B$ is a maximal element of
$\mathcal{U}$ because if $D \in \mathcal{U}$ and 
$B \subseteq D$ then $D \in \mathcal{U}^{\prime}$ and so 
$B = D$. 

\item\label{rem_main_it2}
There is a unital C*-algebra $A$ with 
a maximal UHF~subalgebra such that $MU(A)$ is not 
a maximal element of $(\mathcal{U}, \subseteq)$ 
as in \eqref{rem_main_it1}. For instance, consider 
the separable, simple, unital C*-algebra 
$A = \mathrm{M}_{2^{\infty}} \rtimes_{\alpha} \mathbb{Z}_{2}$ 
where 
$\alpha : \mathbb{Z}_{2} \curvearrowright \mathrm{M}_{2^{\infty}}$ 
is an action with the Rokhlin property. Then by 
\cite[Theorem~3.5]{iz04}, $\mathrm{M}_{2^{\infty}}$ is 
a unital subalgebra of $A$, 
$\mathrm{M}_{2^{\infty}} \cong A$, and 
$\mathrm{M}_{2^{\infty}} \neq A$. Therefore 
$MU(A) = \mathrm{M}_{2^{\infty}}$ is 
a maximal UHF~subalgebra of \( A \), however, 
the only maximal element of \( \mathcal{U} \) 
as in \eqref{rem_main_it1} is \( A \) itself.

\item\label{rem_main_it3}
All maximal elements of $\mathcal{U}$ 
as in \eqref{rem_main_it1} may not be isomorphic. 
For instance, the universal unital free product 
$A = \mathrm{M}_2 \ast_{r} \mathrm{M}_3$ has at least 
two nonisomorphic maximal UHF~subalgebras. For this, put
\[ \mathcal{U}_{\mathrm{M}_{2}} = 
\left\{ D \in \mathcal{U} : \mathrm{M}_2 \subseteq D \right\}
\ \text{and} \ \ \mathcal{U}_{\mathrm{M}_{3}} = 
\{ D \in \mathcal{U} : \mathrm{M}_3 \subseteq D \}. \]
Let $B_{1}$ and $B_{2}$ be maximal elements of 
$(\mathcal{U}_{\mathrm{M}_{2}}, \subseteq)$ and 
$(\mathcal{U}_{\mathrm{M}_{3}}, \subseteq)$, respectively. If
$B_{1} \cong B_{2}$ then since $\mathrm{M}_{2}$ and  
$\mathrm{M}_{3}$ embed into $B_{1}$, $\mathrm{M}_{6}$ 
embeds into $B_{1}$ and hence into $A$, which is a contradiction 
(see Example~\ref{exa_mu_not}). Therefore, 
$B_{1} \ncong B_{2}$. Note that \( B_{1} \) and \( B_{2} \) are 
maximal elements of \( \mathcal{U} \). 
\end{enumerate}
\end{remark}

\textbf{Acknowledgment.} We thank Mikael R{\o}rdam and
David Handelman for very helpful comments and discussions
via email.


\end{document}